\documentclass[runningheads]{llncs}
\usepackage{latexsym}
\def\ra{\rightarrow}

\def\be{\begin{equation}}
\def\ee{\end{equation}}
\def\bea{\begin{eqnarray}}
\def\eea{\end{eqnarray}}
\def\eps{\epsilon}
\def\F{{\bf F}}
\def\Z{{\bf Z}}
\def\C{{\bf C}}
\def\CC{{\cal C}}
\def\MM{{\cal M}}
\def\N{{\bf N}}
\def\Q{{\bf Q}}
\def\R{{\bf R}}
\def\Pr{{\bf P}}
\def\X{{\rm X}}
\def\Aut{{\rm Aut}}
\def\PSL{{\rm PSL}}

\def\SL{{\rm SL}}
\def\GL{{\rm GL}}
\def\Sym{{\rm Sym}}
\def\LogN{\log^{O(1)}\!N}
\def\LogkN{\log^{O_k(1)}\!N}
\def\LogMN{\log^{O_M(1)}\!N}
\def\LogX{\log^{O(1)}\!X}
\def\LogH{\log^{O(1)}\!H}
\def\MW{Mordell-\mbox{\kern-.16em}Weil}
\begin{document}

\title{Rational points near curves and small nonzero~$|x^3-y^2|$
via lattice reduction}
\titlerunning{Points near curves and Hall's conjecture via
lattice reduction}
\author{Noam D.~Elkies}
\institute{Department of Mathematics, Harvard University,
Cambridge, MA 02138 USA\\{\tt elkies@math.harvard.edu}}
\maketitle

\begin{abstract}
We give a new algorithm using linear approximation and lattice
reduction to efficiently calculate all rational points of small height
near a given plane curve~$C$.  For instance, when $C$\/ is the Fermat
cubic, we find all integer solutions of $|x^3+y^3-z^3|<M$\/ with
$0<x\leq y<z<N$ in heuristic time $\ll(\LogN) M$\/ provided
$M\gg N$, using only $O(\log N)$ space.  Since the number of solutions
should be asymptotically proportional to $M \log N$\/ (as long as
$M<N^3$), the computational costs are essentially as low as possible.
Moreover the algorithm readily parallelizes.
It not only yields new numerical examples but leads to
theoretical results, difficult open questions, and natural
generalizations.  We also adapt our algorithm to investigate Hall's
conjecture: we find all integer solutions of $0<|x^3-y^2|\ll x^{1/2}$
with $x<X$\/ in time $O(X^{1/2}\LogX)$.  By implementing this
algorithm with $X=10^{18}$ we shattered the previous record for
$x^{1/2}/|x^3-y^2|$.  The $O(X^{1/2}\LogX)$ bound is rigorous; 
its proof also yields new estimates on the distribution mod~$1$
of $(cx)^{3/2}$ for any positive rational~$c$.
\end{abstract}
\section{Introduction}
One intriguing class of Diophantine problem concerns small values
of homogeneous polynomials.  In the simplest nontrivial case of
a polynomial in three variables defining a projective plane
curve~$C: P(X,Y,Z)=0$, the problem can be reformulated thus:
given a plane curve~$C$, describe for each positive $N,\delta$
the rational points of height at most~$N$\/ in~$\Pr^2$ which are
at distance at most $\delta$ of~$C$.  With present-day methods,
hardly any nontrivial results can be proved on the number
or existence of such points.  But one can still seek numerical
evidence, and efficient algorithms for obtaining this evidence.  
The direct approach is to try all $x,y$ with $|x|,|y|\leq N$,
and for each pair to solve $P(x,y,z)=0$ for~$z\in[-N,N]$, recording
those cases in which $z$ is sufficiently close to an integer.
This requires space $O(\log(N))$ but time $(N^2+\delta N^3)\LogN$,
which is inefficient once $\delta$ is much smaller than $N^{-1}$
since for general~$C,N,\delta\gg N^{-3}$ the number of solutions
should be proportional to $\delta N^3$.  We give a new
algorithm, also requiring only $O(\log(N))$ space, but with
heuristic running time $(N+\delta N^3)\LogN$.  Thus as long
as $\delta\gg N^{-2}$ we expect to find all the points of height
$\leq N$\/ and distance $\leq\delta$ in time only $\LogN$\/
per point.  Moreover, our method readily parallelizes, since it
divides the computation into many independent subproblems.

We describe this algorithm, give the heuristic estimate for its
run time, and briefly discuss the problem, which seems
quite difficult, of proving our heuristic time estimates.
We prove (Thm.\ref{connect_the_dots}) that an alternative
description of those points can always be computed in the
heuristically expected time.  We then discuss natural
generalizations to other valuations and higher dimensions.

An algorithm for finding rational points {\em near} a variety can
in particular find rational points {\em on} the variety; applying
our methods to embeddings of the variety in projective spaces
of high dimension we obtain a new approach to this fundamental
problem in computational number theory which improves on existing
methods in several important cases.  This approach also works
for non-algebraic varieties, and even yields a theoretical
result (Thm.\ref{transcendental}) on the paucity of rational points
on non-algebraic analytic curves.

We next describe experimental results of the implementation
of our algorithm to various curves of interest, notably
the Fermat curves of degree $n>2$, where some of our experimental
findings led us to new polynomial families of small values of
$|z^n-y^n-x^n|$ (Thm.\ref{fermat}).  We devote a separate section
to the case of the cubic Fermat curve, corresponding to small values
of $|z^3-y^3-x^3|$, a problem for which there is already some
literature and the heuristics are subtler.  In particular, we found
for several integers~$d<10^3$ the first representation of~$d$\/ as
a sum of three integer cubes; D.J.Bernstein has since extended the
search up to $N=2\cdot 10^9$ and beyond, and found many new solutions,
including one for $d=30$ which was a long-standing open problem.

Finally we show how to modify our algorithm to efficiently search for
small nonzero values of $|x^3-y^2|$.  This is the topic of Hall's
conjecture, which is part of a web of important Diophantine problems
surrounding the ABC conjecture of Masser and Oesterl\'e.  The
conjecture asserts that $x^3-y^2$ is either zero or
$\gg_\eps x^{1/2-\eps}$ for all $x,y\in\Z$.  We are able
to find all solutions of $0<|x^3-y^2|\ll x^{1/2}$ with $x\leq X$\/
in time $O(X^{1/2}\LogX)$, again using only $O(\log X)$
space.  Using this improvement on the obvious $X\LogX$ method
of trying all $x\leq X$, we computed all cases of $0<|x^3-y^2|<x^{1/2}$
with $X\leq 10^{18}$.  We found ten new solutions, including
most notably
$$
5853886516781223^3 - 447884928428402042307918^2 = 1641843
$$
with $x^{1/2}\big/|x^3-y^2|=46.600+$, improving the previous record by
a factor of almost~$10$.  In this case the time estimate is
{\em not}\/ heuristic; its proof not only streamlined the computation
but even yields new theorems on the distribution mod~$1$ of
$(cx)^{3/2}$ for any positive rational~$c$.  We announce some of these
results at the end of the present paper; the full statements and proofs
will appear elsewhere.

\subsection{Acknowledgements}
Richard K.~Guy wrote the book~\cite{Guy} that first introduced me
to many open problems in number theory including the Diophantine
equations $x^3+y^3+z^3=d$ \cite[Prob.~D5]{Guy}, and later brought me
up to date on recent work on this problem.
Dan J.~Bernstein efficiently implemented my new algorithm for the problem.
Alan Murray told me of the appearance of approximate integer solutions
of $x^{12}+y^{12}=z^{12}$ on {\em The Simpsons.}
Frits Beukers and Franz Lemmermeyer filled gaps in my knowledge
of earlier work concerning Hall's conjecture.
Barry Mazur suggested that a method for locating points near a variety
might also profitably be applied to finding points on the variety;
this started me thinking in the direction that led to Theorems
\ref{2overM} through \ref{transcendental}.
Alf van der Poorten and Hugh Montgomery directed me to
Bombieri and Pila's work~\cite{BP} concerning integral points on curves;
Peter Sarnak put me in contact with Pila, who noted his more recent
paper~\cite{Pila}; meanwhile Victor Miller alerted me to results
announced by Roger Heath-Brown~\cite{HB:nonarch},
who discussed his and Pila's work with me.
Meanwhile, Michel Waldschmidt informed me of relevant results
by Weierstrass and others collected in \cite[Chapter 3]{Mahler}.
I thank them all for these contributions to the present paper.

Most of the numerical and symbolic computations reported here were 
carried out using the {\sc gp/pari} and {\sc macsyma} packages.

This work was made possible in part by funding from the
David and Lucile Packard Foundation.

\section{The algorithm in theory}

\subsection{Specification and heuristic analysis}
While we are mainly interested in algebraic plane curves~$C$,
the algorithm does not require so strong a hypothesis: we can
find\footnote{
  Our computations indicate that the first example is probably the
  smallest value of $|x^\pi+y^\pi-z^\pi|$ for positive integers $x,y,z$,
  and at any rate the smallest with $z\leq 10^6$; and the second is the
  smallest ratio of $|x^7+y^7-z^7|$ to $z^7$, and even to $z^4$,
  for positive integers satisfying $x\leq y < z \leq 10^6$.  See the
  next section.
  }
$2063^\pi + 8093^\pi - 8128^\pi = 0.019369-$ as well as
$386692^7 + 411413^7 = (1-1.035\ldots\cdot 10^{-18}) 441849^7$.
All we need is that $C$\/ is the image of a differentiable map
$\phi:[0,1]\ra\R\Pr^2$ with bounded second derivatives.
Fix a positive $\delta\leq1$, and assume $\delta\gg N^{-2}$ for
reasons given in the next paragraph.  Partition $[0,1]$ into
$O(\delta^{-1/2})$ intervals $I_m$ each of length
$|I_m|=O(\delta^{1/2})$.  On each $I_m$, approximate $\phi$ to within
$O(|I_m|^2)=O(\delta)$ by a linear approximation $\bar\phi$.
Then a point at distance $\leq\delta$ from $\phi(I_m)$ remains
at distance $\ll\delta$ from $\bar\phi(I_m)$.

We now treat each $I_m$ independently.  The triples
$(x,y,z)\in\Z^3-\{{\bf0}\}$ such that $(x:y:z)\in\Pr^2$
has height $\leq N$\/ and is within $O(\delta)$ of $\phi(I_m)$
are among the nonzero integer points in a parallelepiped $P_m$
of height, length and width proportional to $N,\delta^{1/2}N,\delta N$.
Thus we expect that $|P_m\cap\Z^3|$ is approximately the volume
of~$P_m$, provided that this volume is $\gg 1$.  This is the case
once $\delta\gg N^{-2}$.  (That is why we insisted that
$\delta\gg N^{-2}$: choosing smaller $\delta$
would only make us work at least as hard to find fewer points.)
Listing all the points in $P_m\cap\Z^3$ is a standard application
of lattice reduction.  Let $M_m$ be an invertible $3\times 3$
matrix such that $M_m P_m$ is the cube $K=[-1,1]^3$.  We are then
seeking all $v\in\Z^3$ such that $M_m v \in K$, or equivalently
all vectors in $K\cap M_m^{-1}\Z^3$.  We find them by reducing
the lattice $M_m^{-1}\Z^3$.  This gives us a matrix $L_m\in\GL_3(\Z)$
such that $M_m L_m$ is small.  Now $M_m v \in K$\/ if and only if
$w\in\Z^3\cap(M_m L_m)^{-1}K$\/ where $v=L_m w$.  But $(M_m L_m)^{-1}K$
is contained in the box centered on the origin whose $i$-th side is
twice the $l^1$ norm of the $i$-th row of $(M_m L_m)^{-1}$ ($i=1,2,3$).
For each nonzero integral $w$ in this box, calculate $(x,y,z) = L_m w$
and test whether $(x:y:z)$ in fact has height $\leq N$\/ and lies
within~$\delta$ of~$C$.  Doing this for each~$m$ yields the full list
of such points.

As advertised, the algorithm requires only $O(\log N)$ space
(though much more space is usually needed to store the results
of the computation).  Also, since each of many intervals $I_m$ is
treated independently, the computation can be massively parallelized
with little loss among processors that interact only by reporting
each $(x:y:z)$ to headquarters as it is found.  How long do we
expect the computation to take?  We assume that $\phi$ and its
derivatives can be calculated to within $N^{-O(1)}$ in time
$\LogN$.  Such is the case for all curves we consider and for every
algebraic plane curve.  Then each $M_m$ takes only $\LogN$ operations
to compute.  Each lattice reduction can also be done in time polynomial
in~$\log N$, since our lattices are in fixed dimension --- and moreover
our dimension of~$3$ is small enough that Minkowski reduction is
described explicitly.  [For an overview and further references
concerning Minkowski reduction, see \cite[pp.396--7]{SPLAG}.]
So far this amounts to $\delta^{-1/2}\ll N$\/
time up to the usual log factors.  Now each $P_m$ has volume
$2^3/|\det M_m|\ll(\delta^{1/2}N)^3$.  If each $(M_m L_m)^{-1}$
had all of its entries $O(\delta^{1/2}N)$ --- equivalently,
if the shortest nonzero vectors of each lattice $M_m^{-1}\Z^3$
had length $\gg\delta^{-1/2}/N$ --- then there would only be
$O((\delta^{1/2}N)^3)$ choices for~$w$, which summed over~$m$
gives $O(\delta N^3)$.  Thus the total work would indeed be
$\LogN$\/ times the expected number of solutions.
Unfortunately it is too optimistic to expect that the entries of
$(M_m L_m)^{-1}$ are all $\ll\delta^{1/2}N$.  If the lattices
$M_m^{-1}\Z^3$ are randomly distributed in the space of lattices
of covolume $(\delta^{1/2}/N)^3$ in~$\R^3$, some of them will
have nonzero vectors much shorter than $\delta^{-1/2}/N$.  However,
the {\em average} number of lattice vectors in~$K$\/ of a random
lattice of determinant~$D$\/ is still $O(1/D)$.  Thus we expect ---
and typically find in practice --- that, even accounting for the
occasional short lattice vector, we will find all rational points
of height $\leq N$\/ that lie within~$\delta$ of~$C$, doing on
average $\LogN$\/ work per point.

\subsection{Can the estimates be made rigorous?}  Our assumption
that the lattices $M_m^{-1}\Z^3$ are randomly distributed was
not proved; indeed it is false at least for some choices of~$C$.
Most glaringly, if $C$\/ is a rational straight line then there are
$\gg N^2$ rational points on~$C$, and {\em a fortiori} at least as many
at distance $\leq\delta$.  While we of course will not apply our
algorithm to straight lines, we do apply it to the $n$-th Fermat curve,
which has contact of order~$n$ with several rational lines
such as $y=z$; each of those lines contains $\gg N^{2-2/n}$ points
at distance $\ll 1/N^2$ from the curve, exceeding the expected
count of $N\LogN$ once $n>2$.  (These are the points we exclude by
imposing the inequality $y<z$ in $0<x\leq y<z<N$.)  Assume, then,
that $C$\/ has at most finitely many tangent lines which have
contact of order $>2$ with~$C$, and for any $\delta>0$ let
$C_\delta$ be the curve consisting of points of~$C$\/
at distance $\geq\delta$ from each of those higher-order tangent lines.
For each point~$P$\/ on~$C_\delta$ we obtain a lattice
$L_\delta(P)\subset\R^3$ whose nonzero short vectors correspond
to points near~$P$\/ in~$\Pr^2(\Q)$, of height $\ll\delta^{-1/2}$,
lying at distance $\ll\delta$ from~$C_\delta$.
This gives a map $\Lambda_\delta$ from $C_\delta$
to the moduli space of lattices in~$\R^3$.  We would thus like to ask:
as $\delta\ra 0$, does the image of $\Lambda_\delta$ become uniformly
distributed in this moduli space?

There are several problems with this formulation of our question.
A minor one is that we have not defined $\Lambda_\delta$ precisely
enough for the question to make sense, because we have left some
$O$-constants unspecified.  This did not matter for qualitative
properties such as whether the lattice has $O(1)$ short vectors,
but makes it easy to frustrate uniform distribution by simply choosing
$\Lambda_\delta$ to avoid a small region in the moduli space.  This
problem is easy enough to fix for any given $C$\/; for instance, if
$C$\/ is given by $x \mapsto (x:y(x):1)$ for some differentiable
function $y:[0,1]\ra[-1,1]$ with bounded second derivatives,
we may take for $\Lambda_\delta(x)$ the integer span of the
columns of
\be
\left(
\begin{array}{ccc}
0 & 0 & \delta \\
1 & 0 & -x \\
-y'/\delta & \;\, 1/\delta \,\; & (xy'-y)/\delta
\end{array}
\right)
\; .
\label{Lambda(x)}
\ee
But this brings us to a more serious difficulty.  The question of
whether $\Lambda_\delta(C_\delta)$ is asymptotically uniformly distributed
as $\delta\ra 0$ is likely to be a very hard problem in analytic
number theory.  For our purposes we are only concerned with how
often and how close does $\Lambda_\delta(P)$ come near the cusp of
the moduli space.  For instance, we see in the final section
that if $C$\/ is a conic then $\Lambda_\delta(C_\delta)$ is restricted
to a surface in the moduli space of lattices in~$\R^3$, but within
that surface it still approaches the cusp rarely enough that
the average number of short vectors in a lattice in $\Lambda_\delta(C)$
is still $\ll \log(1/\delta)$.  In general, then, what we would like
is the following result: as $\delta\ra 0$, the average number of vectors
of norm $<1$ of a lattice in $\Lambda_\delta(C_\delta)$ is
$\ll \log^{O(1)}(1/\delta)$.

This still looks like a very difficult problem.  While it remains open,
we propose a contingency plan in case the lattices $L_\delta(P)$
have many more short vectors than expected.
If all the short vectors are multiples of a single vector of small norm,
there is no difficulty, because all these multiples yield the same point
in~$\Pr^2$.  But there could be two independent short vectors, whose
linear combinations yield a line in~$\Pr^2$ containing many points of
small height near~$C$.  We claim that this is in fact the only way
that a lattice of covolume $\ll 1$ could have more than $O(1)$ short
vectors.  This claim is easy enough to check using the description of
Minkowski-reduced lattices in~$\R^3$, but we shall later need a
generalization to lattices in higher dimension.  We thus state
and prove the generalization as follows:

\begin{lemma}
For each positive integer $n$ and positive real $t$ there exists
an effective constant $M_n(t)$ such that the following bound holds:
for any lattice~$\Lambda\subset\R^n$ whose dual lattice $\Lambda^*$
has no nonzero vector of length $<r$, and for any $R>0$, there are
at most $M_n(rR) \, r^{-n} |\Lambda|^{-1}$ vectors of length~$\leq R$ 
in~$\Lambda$.
\label{reduction}
\end{lemma}

Here $|\Lambda|$ is the covolume Vol$(\R^n/\Lambda)$.  The lemma
can be obtained as a consequence of the theory of lattice reduction,
but it is not easy to extract $M_n(t)$ explicitly this way.  We thus
give the following alternative proof in the spirit of~\cite{Cohn} from
which explicit (albeit far from optimal) bounds $M_n(t)$ may be
easily computed if desired.

\begin{proof}
Given $n$, choose a {\em positive} Schwartz function $f:\R^n\ra\R$
with the following properties:
$f$\/ is {\em radial}, i.e.\ $f(x)$ depends only on~$|x|$; 
and the Fourier transform $\hat f: \R^n\ra\R$, defined for $y\in\R^n$ by
\be
\hat f(y) := \int_{x\in\R^n} f(x) e^{2\pi i (x,y)} dx,
\label{Fourier}
\ee
satisfies $\hat f(y)\leq0$ for all $y$ such that $|y|\geq 1$.
For instance, we may take
\be
f(x) = (|x|^2 + a) e^{-\pi c |x|^2}
\label{sample_f}
\ee
where
\be
0 < c < \frac{2\pi}{n}, \quad
a = \frac{1}{c^2} - \frac{n}{2\pi c},
\label{a,c}
\ee
because the Fourier transform of a function (\ref{sample_f}) is
\be
\hat f(y) = \left( a + \frac{n}{2\pi c} - \frac{|y|^2}{c^2} \right)
e^{-\pi |y|^2 / c}
\label{fhat}
\ee
for any $c>0$ and $a\in\R$.  By Poisson summation,
\be
\sum_{x\in\Lambda} f(rx) =
\frac1{r^n|\Lambda|} \sum_{y\in\Lambda^*} \hat f(y/r).
\label{poisson}
\ee
Under the hypothesis on~$r$, the only positive term in the sum
over~$y$ is $\hat f(0)$.  The sum over~$x$ is bounded from below
by the sum over~$x$ of length $\leq R$, which is at least
the number of such vectors times $\min_{|x|\leq R}^{\phantom0} f(rx)$.
It follows that $\Lambda$ has at most
\be
\frac{\hat f(0)} {r^n |\Lambda| \min_{|x|\leq rR} f(x)}
= M_n(rR)\, r^{-n} |\Lambda|^{-1}
\ee
vectors of length $\leq R$, as claimed.
\end{proof}

\begin{corollary}
For each positive integer $n$ there exists an effective constant
$A_n$ such that if a lattice $\Lambda\subset\R^n$ has more than
$A_n R^n / |\Lambda|$ vectors of length $<R$\/ for some $R>0$
then all those vectors lie in a hyperplane, which can be computed
in polynomial time.
\label{reduction_cor}
\end{corollary}

\begin{proof}
Except for the last phrase, this follows from the previous Lemma
by taking $r=1/R$\/ and $A_n=M_n(1)$, since then $\Lambda^*$ must
have a nonzero $y$ of length at most~$r$, and any vector of~$\Lambda$
of length $<R$\/ must be orthogonal to~$y$.  To assure that $y$ can be
computed in polynomial time, we take $r=c/R$\/ for a positive
constant~$c$ small enough that if $\Lambda^*$ has a nonzero vector
of length at most~$c/R$\/ then the LLL algorithm will find a
(possibly different) nonzero vector of length at most $1/R$.
Our Corollary now holds with $A_n = c^{-n} M_n(c)$.
\end{proof}

{}From the case $n=3$ of this Corollary we deduce:

\begin{theorem}
Let $C$\/ be the image of a differentiable map
$\phi:[0,1]\ra\R\Pr^2$ with bounded second derivatives.
Then for each $N>1$ and $\delta\geq N^{-2}$ one can find
$O(\delta N^3)$ rational points and $O(N)$ rational line
segments each of length $O(1/N)$ in $\Pr^2$ which together
include all rational points of height $\leq N$\/ at distance
$\leq \delta$ from~$C$.  These points and line segments
can be computed in time $\ll \delta N^3 \LogN$.  Outside
of $O(\delta N^3 \log N)$ space used only to record each point or
segment as it is found, the computation requires space
$\ll\LogN$.  All implied constants depend effectively on~$C$.
\label{connect_the_dots}
\end{theorem}

Note that here we do not exclude neighborhoods of high-order
rational tangents to~$C$\/; such tangents will contain 
some of the lines segments computed by the algorithm.

\subsection{Variations and generalizations}

The problem of finding rational points near plane curves is only
the first nontrivial example of many analogous problems to which
our method can apply.  We briefly discuss some of these here.

One easy variation is to change the norm: instead of approximating
the curve in the real valuation, use a nonarchimedean one, or a
combination of several.  For instance, one can efficiently seek
nontrivial triples of small integers the sum of whose cubes is divisible
by a high power of~$2$ or of~$10$.  Likewise one can replace $\Z$
by $\F_q[T]$ or similar rings in function fields of positive genus.
The lattice-reduction step should then be even easier
than in the archimedean case, though in the function-field setting
our approach faces strong competition from the method
of undetermined coefficients, and it is not clear which is superior.  
All these comments apply equally to the adaptation of our method
to the problem of finding small nonzero values of $|x^3-y^2|$,
provided the characteristic is not~$2$ or~$3$.  For the $|x^3-y^2|$
problem, the work estimates are again rigorous; otherwise, they are
still heuristic, but their analysis may be more tractable in the
function-field case.

Higher dimensions present many new opportunities.  The easiest
generalization is to a $\CC^2$ hypersurface in $\Pr^{k-1}$.  Here
we are seeking small values of a homogeneous function of $k$ variables
evaluated at an integral point.  This time we chop the hypersurface
into $O(\delta^{-(k-2)/2})$ chunks each of diameter $O(\delta^{1/2})$,
and replace each chunk by a subset of a hyperplane which approximates
it to within $O(\delta)$.  The points of height $\leq N$\/ that are
within $O(\delta)$ of this chunk then come from integral points
in a parallelepiped in~$\R^k$ whose sides have lengths $O(N)$,
$O(N\delta^{1/2})$ ($k-2$ times), and $O(N\delta)$.  Again most of
these this parallelepipeds have volume $\gg 1$ provided
$\delta\gg N^{-2}$, and we locate the integral points
using lattice reduction in $\R^k$.  So, as long
as $\delta\gg N^{-2}$, we expect to find on the order of $\delta N^k$
points, using $\ll\LogkN$ space and spending $\ll\LogkN$ time per
point.  For a general hypersurface, this again improves on other
approaches to the problem.  But the improvement decreases with~$k$\/:
the direct approach takes time $N^{k-1}\LogN$, and we lower the
exponent by a factor no better than $(k-2)/(k-1)$, which approaches~$1$
as $k\rightarrow\infty$.  Moreover, lattice reduction in~$\R^k$
quickly becomes difficult as $k$\/ grows.  Another consideration
is that for special surfaces there are known, and simpler, algorithms
that take time $N^{k-1}\LogN$\/ or less once $k>3$.  For instance,
for Fermat surfaces in $\Pr^3$, one readily adapts the method
of~\cite{Bernstein} to find all solutions of
$x^n+y^n=z^n\pm t^n+O(z^{n-2})$
in positive integers with $t\leq z\leq N$,
in expected time $N^2\LogN$, and with no need for lattice
reduction in~$\R^4$ or other complicated ingredients.
This computation does require space proportional to~$N\log N$,
which however poses no difficulty for practical values of~$N$.
As in the previous paragraph, all that is described in the present
paragraph can be done also for a nonarchimedean norm, with
similar results except that lattice reduction over a function field
is tractable even for large~$k$.  In either case rigorous estimates
may become even less accessible as $k$ grows.

We can generalize further to manifolds $\MM\subset\Pr^{k-1}$
of codimension $c>1$.  Here we expect to find on the order of
$\delta^c N^k$\/ rational points of height $\leq N$\/
at distance $O(\delta)$ from~$\MM$.  We chop $\MM$\/ into
$O(\delta^{-(k-1-c)/2})$ patches of diameter $O(\delta^{1/2})$,
each of which yields a parallelepiped in~$\R^k$ with dimensions
of order~$N$\/ (once), $N\delta$ ($d$\/ times), and $N\delta^{1/2}$
(the remaining $k-1-c$ dimensions).  We thus expect to efficiently find
all $\sim \delta^c N^k$ points as long as $\delta\gg N^{-2k/(k+c-1)}$.
A further possibility emerges if $\MM$\/ has bounded derivatives
past the second derivatives and has small enough dimension
compared with~$k$\/: we can then make further headway when $\delta$
falls below that threshold.  Usually we are only interested in points
much closer than $N^{-2k/(k+c-1)}$; but as long as we use only the
$\CC^2$ structure we gain nothing by making $\delta$ even smaller,
so we may as well find all the points at distance $O(N^{-2k/(k+c-1)})$
and locate the best approximations in the resulting list.
However, if $\MM$\/ is $\CC^3$ and its
dimension $d=k-1-c$ is so small that $k>{{d+2}\choose 2}$, then
a patch of diameter $\eps$ is contained in a box with
$d$\/ sides of length $\ll\eps$, a further $(d^2+d)/2$ sides
of length $\ll\eps^2$, and the remaining $k-{{d+2}\choose 2}$
sides of length $\ll\eps^3$.  This means that we can makes our
parallelepipeds thinner in some directions, and thus use wider
patches of~$\MM$, covering the entire manifold with fewer of them.
This lets us locate the points of height $\leq N$\/ closest to~$\MM$\/
in time significantly less than it would take to record all
the points at distance $\ll N^{-2k/(k+c-1)}$, even though not
so efficiently that we only spend $\ll\LogkN$ time per point.
More generally if $\MM$ is a $\CC^i$ manifold we can exploit
bounds on the $i$-th derivatives once $k>{{d+i-1}\choose{i-1}}$.

If the ambient projective space is not of high enough dimension,
we can still make some use of approximations to~$\MM$\/ of degree
$i>1$ by using the $i$-th Veronese embedding~$V_i$ of $\Pr^{k-1}$ into
projective space of dimension ${{k+i-1}\choose{i}} - 1$.  [The $i$-th
{\em Veronese embedding} takes the point with projective coordinates
$(X_1:\cdots:X_k)$ to the point whose projective coordinates are
all ${k+i-1}\choose{i}$ monomials of degree~$i$ in the $X_j$.
Thus $V_i$ raises all heights to the power~$i$, and transforms
intersections with hypersurfaces of degree~$d$ in~$\Pr^{k-1}$
into hyperplane sections in a projective space of much higher
dimension.  For more on Veronese embeddings, see for instance
\cite{FH}, where they arise several times.]
The idea is to surround each patch of $V_i(\MM)$ by a box containing
all points in $V_i(\Pr^{k-1})$ at distance $O(\delta)$ from $V_i(\MM)$.
The resulting asymptotic improvement may be only barely worth it
in practice, though.  Consider the simplest case of a $\CC^3$ curve
$C\in\Pr^2$, embedded in $\Pr^5$ by $V_2$.  Assume for simplicity
that the parametrization $\phi$ of~$C$\/ has $|\phi''|$ bounded
away from zero.  Then, for $\delta$ such that
$\eps^3\ll\delta\ll\eps^2$, the radius-$\delta$ neighborhood in~$\Pr^2$
of an interval of length~$\eps$ on~$C$\/ maps into a box in~$\Pr^5$
whose sides are of order $\eps,\eps^2,\delta,\eps^3,\eps^4$.
[To see this, choose coordinates $(X_0:X_1:X_2)$ on~$\Pr^2$ for which
$\phi$ is of the form $(1:t:t^2+O(t^3))$ for~$t$ in a neighborhood
of~$0$, and note that $V_2$ takes $(X_0:X_1:X_2)$ to
$(X_0^2 : X_0 X_1 : X_0 X_2 : X_1^2 : X_1 X_2 : X_2^2)$.]
Thus the points of height at most~$N$\/ in that neighborhood map
to lattice points in a $6$-dimensional parallelepiped of volume
$\ll \delta N^{12}\epsilon^{10}$.  (Here $N$\/ occurs to the power $12$
rather than $6$ because $V_2$ squares the height of each rational
point.)  Thus if we take $\eps=(\delta N^{12})^{-1/10}$
we expect to find all points at distance $\ll\delta$ from~$C$,
of which there should be about $\delta N^3$, in time
$N(\delta N^2)^{1/10}\LogN$.  The condition $\delta\gg\eps^3$ yields
$\delta\gg N^{-36/13}$, so we save a factor of at most $N^{1/13}$.
We pay not only by missing the points at distance between $N^{-36/13}$
and $N^{-2}$ (which usually do not interest us anyway) but also by
reducing lattices of rank~$6$ rather than~$3$.  This takes more time
per lattice, and probably yields parallelepipeds whose average bounding
box is larger.  Each of these effects amounts to only a constant factor,
but these factors may be considerable, and it will be interesting to see
how large $N$\/ must be for this use of $V_2$ to be practical.

\subsection{Rational points on varieties}

In the last paragraph we exploited the fact that points near~$\MM$\/
map under $V_i$ to points that are not only near $V_i(\MM)$ but exactly
on $V_i(\Pr^{k-1})$.  We can go much further when we search for points
exactly on~$\MM$.  Again we consider the simplest case of a curve.
We begin with a curve in one projective space:

\begin{theorem}
Let $C$ be an algebraic curve in $M$-dimensional projective space,
defined over~$\Q$ and not contained in any hyperplane.
Then for any $N\geq1$ the rational points of~$C$ of height at most~$N$\/
can be listed in time $\ll_C N^{2/M} \LogMN$.  The implied constants
depend effectively on~$d$ and~$C$.
\label{2overM}
\end{theorem}

\noindent
{\em Remarks.}  As seen above for $M=2$, this result applies more
generally to a $\CC^M$ curve in $\Pr^M$ whose intersection with
any hyperplane can be computed in polynomial time.
The exponent $2/M$\/ is best possible: a rational normal curve
of degree~$M$\/ (a.k.a.\ the image of $\Pr^1$ under $V_M$)
has on the order of $N^{2/M}$ rational points of height at most~$N$,
and it takes time $\gg N^{2/d} \log N$\/ just to write them down.
The constant implied in $O_d(1)$ and/or $\ll_C$, while effective,
may be unpleasant in practice for large~$M$, since
lattice reduction in dimension $M+1$ is involved.

\begin{proof}
A segment of~$C$\/ of length $\ll N^{-2/M}$ is contained in a box
whose $i$-th side is $\ll N^{-2i/M}$ ($i=1,2,\ldots,M$\/).  The
rational points of height at most~$N$\/ in this box come from points
of $\Z^{M+1}$ contained in a box~$B$\/ whose $i$-th side is
$\ll N^{1-2i/M}$ ($i=0,1,2,\ldots,M$\/) and thus has volume $O(1)$.
It takes time $\ll_C \LogMN$\/ to apply lattice reduction and,
by Corollary~\ref{reduction_cor}, either list $\Z^{M+1}\cap B$\/
or find a hyperplane containing $\Z^{M+1}\cap B$.
In the former case, we test whether each of
the resulting $O_d(1)$ points lies in~$C$.  In the latter case,
we map this hyperplane to $\Pr^M$ and intersect it with~$C$,
finding at most $\deg(C)= O_C(1)$ rational points.
Thus in either case we find all rational points of height $\leq N$\/
on our segment in time $\ll_C \LogMN$.  Since it takes only
$O(N^{2/M})$ segments to cover~$C$, we are done.
\end{proof}

It might seem that this algorithm is superfluous: if $C$\/ has
genus~$0$ then its small rational points may be found directly from a
rational parametrization, without any lattice reduction; and if $C$\/
has positive genus then we can find all its points of height $\leq N$\/
in time $\ll \LogN$\/ once we have generators of the \MW\ group of the
Jacobian of~$C$.  But the difficulty is that we must first find
these generators, and this requires locating rational points on
a curve or a higher-dimensional variety.  For instance, to find
the \MW\ group of an elliptic curve~$E$\/ we usually apply a few
descents and then search for points on certain principal homogeneous
spaces for~$E$, each of which is a curve~$C$\/ of genus~$1$, usually
(in the case of a complete $2$-descent) of the form $y^2=P(x)$
for some irreducible quartic $P\in\Z[X]$.  One then searches
for~$x\in\Q$\/ of height up to~$H$\/ for which $P(x)\in\Q^2$.
There are on the order of $H^2$ candidates for~$x$; one can set
up a sieve to efficiently try them all, but this still takes
time $H^2 \LogH$\/ (and significant space).  Instead we can
embed $C$\/ in~$\Pr^3$ as the intersection of two quadrics
(by writing $P(x)$ as a homogeneous quadric in $1,x,x^2$),
and use the algorithm of Thm.\ref{2overM} with $N=H^2$
to find all rational solutions of $y^2=P(x)$ with $x$ of height
$\leq H$\/ in time $H^{4/3} \LogH$.  For certain~$E$\/ one can
use Heegner points to locate a rational point on~$C$\/ to within
$\delta$ (see \cite{NDE:Heegner}); if $\delta\ll H^{-2}$,
this is sufficient to identify $x$ using continued fractions,
a.k.a.~lattice reduction in dimension~$2$.  Using the new
algorithm, we see that $\delta\ll H^{-4/3}$ suffices if we use
lattice reduction in dimension~$4$.  This saves a constant factor
in the computation of~$x$, since fewer digits and terms
are needed in the floating-point computation of Heegner points.
When $C$\/ has genus $>1$, there are only finitely many rational
points by Faltings' theorem, but they still may be of significant
number and/or height.  For instance, in \cite{KK,Stahlke} one finds
curves $C: y^2=P(x)$ of genus~$2$ which have hundreds of rational
points.  In both cases, all points with $x$\/ of height $\leq 10^6$
were found using the $H^2 \LogH$\/ sieve method, a substantial
computation.  At least in the case considered in~\cite{Stahlke},
where the Jacobian of~$C$\/ is absolutely simple with large \MW\ rank,
it would probably be even more onerous to find all these points by
first determining the \MW\ group.  But the embedding
$(1:x:x^2:x^3:y)$ of~$C$\/ into $\Pr^4$ yields an improvement
from $H^2$ to $H^{3/2}$ with $5$-dimensional lattice reduction. 

We can do even better by mapping the same curve to larger projective
spaces.  Fix an algebraic curve~$C$\/ of genus~$g$\/ defined
over~$\Q$, and a divisor $D$\/ on~$C$\/ of degree $d>0$.
For $n$ sufficiently large, the sections of $nD$\/ embed~$C$\/ into
$\Pr^{nd-g}$.  This embedding sends any rational point on~$C$\/
of height (exponential, as usual here) $\leq H$\/ relative to~$D$\/
to a point on $\Pr^{nd-g}$ of height $\ll H^n$.  By Thm.\ref{2overM}
again, we can find all such points in time $\ll H^{2n/(nd-g)} \LogH$.
Letting $n\ra\infty$, we conclude:

\begin{theorem}
Fix an algebraic curve~$C/\Q$\/ and a divisor $D$\/ on~$C$\/ of
degree $d>0$.  For each $\epsilon>0$ there exists an effectively
computable constant $A_\epsilon$ such that for any $H\geq 1$
one can find all points of~$C$\/ whose height relative to~$D$\/
is at most~$H$\/ in time $A_\epsilon H^{(2/d)+\epsilon}$.
\label{2+eps}
\end{theorem}

For instance, all rational points on $y^2=P(x)$ with $x$\/ of height
at most~$H$\/ can be computed in time $\ll_\epsilon H^{1+\epsilon}$.

What of varieties~$\MM$\/ of dimension $\Delta>1$ in $\Pr^M$?
A chunk of radius~$\delta$ then yields the intersection of $\Z^{M+1}$
with a box with sides as follows: one of length~$O(N)$,
$\Delta$ sides of length $O(N\delta)$,
${\Delta+1}\choose 2$ sides of length $O(N\delta^2)$, \ldots,
${\Delta+i-1} \choose i$ sides of length $O(N\delta^i)$, \ldots\
until ${{\Delta+j} \choose j} = \sum_{j=0}^i {{\Delta+i-1} \choose i}$
first exceeds $M$.  As usual we choose $\delta$ so that the product
of these sides is~$1$, and apply lattice reduction to each of
$O(\delta^{-\Delta})$ chunks.  The difficulty here is that if the
lattice is nearly degenerate, the hyperplane found in
Corollary~\ref{reduction_cor} meets $\MM$ not in a finite
number of points but in a subvariety of positive dimension~$\Delta-1$.
This suggests an induction on~$\Delta$, since we can apply our method
to that hyperplane section of~$\MM$.  But already for $\Delta=2$
such an argument requires a version of Thm.\ref{2overM}
with more uniformity in the implied constants than we know how to
obtain.  However, as with our first nontrivial case of curves
in~$\Pr^2$, we do not expect such degenerate lattices to arise
in practice often enough to raise the computational cost above
$O(\delta^{-\Delta} \LogN)$, except for a finite number of proper
subvarieties of~$\MM$.  If we assume this, we can again obtain
better estimates by embedding $\MM$\/ in larger projective spaces.
Fix an ample divisor~$D$\/ on~$\MM$, and ask for all rational
points whose height relative to~$\MM$\/ is at most~$H$.  Using
the sections of $nD$\/ to embed $\MM$\/ in projective spaces,
and letting $n\ra\infty$, we find the following heuristic
generalization of Thm.\ref{2+eps}: for each $\epsilon>0$,
there exists a proper subvariety $\MM_0(\epsilon)$ of~$\MM$\/
such that all points of $\MM - \MM_0(\epsilon)$ of height at most~$H$\/
relative to~$D$\/ can be found in time
\be
O_\epsilon(H^{((\Delta+1)/|D|) + \epsilon}),
\label{Delta+1+eps}
\ee
where $|D|$ is the $\Delta$-th root of the intersection number
$D^\Delta$.  One might even hope that $\MM_0(\epsilon)$ can be taken
independent of~$\epsilon$.  For instance, if $\MM$ is a surface of
degree~$d$\/ in $\Pr^3$ then we expect that, for some union
$\MM_0$ of curves on~$\MM$, we can find all rational points
of height~$\leq N$\/ on $\MM-\MM_0$ in time
$\ll_\epsilon N^{(3/\sqrt d\,) + \epsilon}$.
We must admit that this is unlikely to yield a practical improvement
over the $N^2\LogN$ method we already knew: the first $V_i$
that reduces the exponent of~$N$\/ below~$2$ is $V_3$, and then
(assuming $d\geq 4$) the exponent drops only to $24/13$ ---
but instead of reducing $4$-dimensional lattices we are then faced
with lattice reduction in dimension~$20$.  It will probably be a long
time before $N$\/ can feasibly be taken large enough that this extra
effort is worth the $N^{2/13}$ factor gained.

Returning to plane curves, we can use this idea to prove an even
stronger bound on rational points on a plane curve~$C$\/ that is
analytic but not algebraic.  This is because the homogeneous monomials
of degree~$i$ in the coordinates of~$C$\/ are linearly independent
for each~$i$, so $V_i(C)$ spans a projective space whose dimension
grows quadratically in~$i$ (whereas for an algebraic curve the
growth is always linear).  This leads us to the following result:

\begin{theorem}
Let $C$ be a transcendental analytic arc in\/~$\Pr^2$, i.e.\
$C = \{f(x): a\leq x \leq b\}$ where $f$ is an analytic map from
a neighborhood of~$[a,b]$ to~$\Pr^2$ whose image is contained
in no algebraic curve.  Then for each $\epsilon>0$ there exists
a constant $A_\epsilon$ such that for every $H\leq 1$
there are fewer than $A_\epsilon H^\epsilon$ points of height~$\leq H$
in $C \cap \Pr^2(\Q)$.
\label{transcendental}
\end{theorem}

\begin{proof}
For each positive integer~$i$ consider $V_i(C)\subset\Pr^{(i^2+3i)/2}$.
Since $C$\/ is transcendental, $V_i(C)$ is an analytic arc
$V_i\circ f$\/ contained in no hyperplane of $\Pr^{(i^2+3i)/2}$.
Now apply the argument for Thm.\ref{2overM} with $N=H^i$.  As noted
in the remarks following the statement of that theorem, the curve
need not be algebraic as long as it is $\CC^M$ and its intersection
with any hyperplane is of bounded size.  (Here we need not compute this
intersection numerically, since we are only bounding the number of
rational points of small height on~$C$, not computing them efficiently.)
The differentiability is clear since $V_i(C)$ is analytic, and the
boundedness is proved in the next lemma.  We conclude that the number
of points of height~$\leq H$\/ on~$C$\/ is $\ll_i H^{4/(i+3)}$.  Since
$i$ can be taken arbitrarily large, our theorem follows.
\end{proof}

The existence of an upper bound on the size of the intersection of
any hyperplane with $V_i(C)$ is a special case of the following lemma
in complex analysis.  Throughout the lemma and its proof we count zeros
of an analytic function according to multiplicity, even though in
the application to Thm.\ref{transcendental} a multiple zero is no
worse than a simple one.

\begin{lemma}
Let $E$ be an open subset of~$\C$ and $V$ a finite-dimensional vector
space of analytic functions: $E\ra\C$.  Then for any compact subset
$K\subset E$ there exists an integer~$n$ such that any nonzero $f\in V$
has at most $n$ zeros in~$K$.
\label{degree}
\end{lemma}

\begin{proof}
Fix $K$.  We shall say that a compact $K'\subset E$\/ is ``good''
if its boundary $\partial(K')$ is rectifiable and its interior
$K' - \partial(K')$ contains~$K$.  Choose a good $K_1$, and define
a norm on~$V$\/ by $\|f\| = \sup_{z\in K'} |f(z)|$.  Let $V_1$ be
the unit ball $\{ f\in V : \|f\| = 1 \}$.  It is sufficient to
prove the lemma for $f\in V_1$.

For each $f\in V_1$ choose a good $K_f\subseteq K'$ such that
$f$\/ does not vanish on $\partial(K_f)$.  Let
$r_f^{\phantom0} = \inf_{z\in\partial(K_f)}|f(z)|$, and let $n_f$ be
the number of zeros of~$f$\/ in $K_f$.  By Rouch\'e's theorem,
if $g\in V$\/ with $\|f-g\|<r_f$ then $g$ has at most $n_f$ zeros 
in~$K_f$, and {\em a fortiori} in~$K$.  Now $V_1$ is compact and
is covered by the open balls $B_f$ of radius~$r_f$ about $f\in V_1$.
Thus there is a finite subcover $\{B_{f_i}\}_{i=1}^M$.  Then
$n := \max_i n_{f_i}$ is an upper bound for the number of zeros
in~$K$\/ of any $f\in V_1$, and thus of any nonzero $F\in V$.
\end{proof}

To recover our result on hyperplane sections of $V_i(C)$, take
$K=[a,b]$, let $E$\/ be a neighborhood of~$K$\/ on which
$f$\/ is analytic, and choose any analytic functions
$f_0,f_1,f_2$ on~$E$\/ such that $f=(f_0:f_1:f_2)$ on~$E$.
Then take for $V$\/ the space of homogeneous polynomials of
degree~$i$ in $f_0,f_1,f_2$.  If we understand $f$\/ well enough
to obtain for each~$i$ an effective bound~$n$ in Lemma~\ref{degree}
then the constants $A_\epsilon$ in Thm.\ref{transcendental}
are effective too.

With a little additional work $\Q$ can be replaced by an arbitrary
number field~$F$\/ embedded in~$\C$, and $C$\/ by $f(K)$,
where $K\subset\C$ is any compact subset and $f$\/ is again
an analytic map from a neighborhood of~$K$\/ to~$\Pr^2$
whose image is contained in no algebraic curve.

A separate approach to bounding the number of rational points on curves
was initiated in~\cite{BP} and pursued further in~\cite{Pila} and
\cite{HB:nonarch}.  For example, Heath-Brown obtains
in~\cite{HB:nonarch} bounds on the number of rational points
on an algebraic plane curve that coincide with the time estimates
in our Theorems \ref{2overM} and~\ref{2+eps}.  Moreover,
our Thm.\ref{transcendental} is contained in \cite[Thm.8]{Pila},
which asserts that for a number field~$F$\/ with $[F:\Q]=n$ the
number of $F$\/-rational points of height $<H$\/
on a transcendental analytic arc~$C$\/ is at most
$A_{C,n,\epsilon} H^\epsilon$.
Probably the methods of \cite{BP,Pila} can also prove these results
with arcs~$C$\/ replaced by compact transcendental curves $f(K)$,
and our bounds can also be made uniform in~$F$\/ given $[F:\Q]$.
There is clearly some overlap between the two approaches; for instance
the Corollary preceding \cite[Thm.8]{Pila} is the same as our
Lemma~\ref{degree}, but proved using the determinants of~\cite{BP,Pila}.
What is not clear, but intriguing, is whether those determinantal
methods and our lattice-reduction technique can ultimately
be interpreted as facets of the same basic idea.

All this also suggests the question of whether a transcendental arc
can contain infinitely many rational points, of whatever height.
I thank Michel Waldschmidt for pointing out that this question
was already asked, and later answered affirmatively, by Weierstrass.
See \cite[Chapter 3]{Mahler} for this and related results.

\section{The algorithm in practice}

In this section we report on the outcome of the application of
our algorithm to various plane curves, and on some results
suggested by our findings.  We suppress details of the explicit
constants replacing each $O(\cdots)$ and $\ll$; these details
are of course crucial in practice, but are straightforward
and not enlightening.  In each case our curve has some rational
points of inflection, and we make sure to truncate our curve
enough to avoid the tangents at such points but not so much
that we lose approximations near but not on those tangents.

In general, for a plane curve given by a homogeneous equation
$P(X,Y,Z)=0$ of degree~$n$, we associate to a rational point
$(x:y:z)$ near but not on the curve the number
\be
n \max(|x|,|y|,|z|)^{n-3} \bigl/ |P(x,y,z)| \bigr. ,
\label{merit}
\ee
which measures how close the point $(x:y:z)$ is to the curve relative
to the point's height.  We insert the factor $n$ so that we can
reasonably compare approximations for curves of different degrees.
For instance, for the Fermat curve one expects that as $x,y$ vary,
the integer $z^n := x^n+y^n$ comes on average within $\frac14 n z^{n-1}$
of the nearest $n$-th power of an integer, and thus that the smallest
value of $|z^n-y^n-x^n|$ for $z\in[N,2N]$ is proportional to
$n z^{n-3}$.  One could insert further factors to correct for
the length and shape of our curve, but these factors are not
significant for most of the curves we study.

We noted already that the heuristics leading to formulas such as
(\ref{merit}) refer to ``random'' $(x:y:z)$ near the curve,
not for systematic families of approximations which may attain
values of the ratio (\ref{merit}) larger or more often than
expected.  We again give an example for the Fermat curves,
which were the subjects of most of our computations.  One
usually guesses that for each~$r$ there will be $\ll r\log N$\/
triples $(x,y,z)$ such that the ratio (\ref{merit}) exceeds~$r$.
However, in the identity
\be
(t+1)^n - (t-1)^n = 2n t^{n-1} + O(t^{n-3}),
\label{fermdiff}
\ee
we can make $2n t^{n-1}$ an $n$-th power by setting $t=2n u^n$;
this yields $\gg N^{1/n}$ triples with (\ref{merit}) bounded
away from zero.  We note the special cases $n=2,3$ of this identity:
for $n=2$, the $O(t^{n-3})$ error vanishes, and we recover a familiar
parametrization of Pythagorean triples; for $n=3$, the error is
constant, and we can scale the identity to obtain the known family
of solutions $(x,y,z)=(6t^2,6t^3-1,6t^3+1)$ of $z^3-y^3-x^3=2$.
Returning to general $n$: in our searches we set the threshold on
$z^{n-3}/|z^n-y^n-x^n|$ low enough to find all the examples coming
from~(\ref{fermdiff}), as a check on the computation; but we chose
a higher threshold for the tabulation of results so that our list
is not dominated by this polynomial family.

\subsection{Fermat curves of degree $>3$}

We implemented our algorithm to find small values of
$|z^n-y^n-x^n|$ with $0<x\leq y<z$, $4\leq n\leq 20$,
and $z\in[10^3,10^6]$.  Since the threshold for ``small''
depends on the size of~$z$, we wrote $[10^3,10^6]$ as the
union of $10$ intervals $[N/2,N]$ and treated each separately.
We also used a direct search for $z<5000$, using the overlap region
$[1000,5000]$ as a check on the computation.  We did not attempt
to fine-tune the algorithm for efficiency, since we carried it out
more as a demonstration project than a major computational undertaking.
Thus we programmed the search in~{\tt gp}, using the built-in
arithmetic and LLL lattice reduction.  We estimate that transcribing
the program to~C, and replacing LLL by Minkowski reduction in~$\R^3$,
would speed the computation by roughly an order of magnitude; of
course a machine faster than a Sun Sparcstation Ultra~$1$ would
help too.  With a C program and a more powerful machine, 
it should be feasible to search the range $n\in[4,20]$, $z<10^9$
in time on the order of a month.

The behavior of the run times and the counts of solutions with
$|z^n-y^n-x^n|\ll z^{n-2}$ seem broadly consistent with our heuristics,
though we have not attempted a detailed statistical analysis.  We
tabulate the most striking examples, those with
\be
r := nz^{n-3} \bigl/ (z^n-y^n-x^n) \bigr.
\label{r.def}
\ee
of absolute value at least~$4$:

\vspace*{1ex}

\setlength{\tabcolsep}{5pt}
\centerline{
\begin{tabular}{r|c|c|c|r}
$n$ &    $x$ &   $y$ &   $z$ & $r$\ \
\\ \hline
 4  &    167 &    192 &    215 & $ -4.5$ \\
 4  &   8191 &  16253 &  16509 & $ 12.9$ \\
 4  &  24576 &  48767 &  49535 & $-64.5$ \\
 4  &  49152 &  97534 &  99070 & $- 8.1$ \\
 4  &  34231 & 157972 & 158059 & $  5.2$ \\
 4  &  76215 & 311390 & 311669 & $-14.8$ \\
 5  &     13 &     16 &     17 &$-120.4$ \\
 5  &     26 &     32 &     34 & $-15.1$ \\
 5  &     39 &     48 &     51 & $- 4.5$ \\
 5  &     42 &     71 &     72 & $- 8.8$ \\
 5  &    262 &    328 &    347 & $- 6.2$ \\
 5  &   1125 &   2335 &   2347 & $- 5.0$ \\
 5  &   5088 &  16155 &  16165 & $  4.1$ \\
 5  & 190512 & 292329 & 298900 & $  5.5$ \\
 6  &   1236 &   3587 &   3588 & $ 12.5$ \\
 6  &   6107 &   8919 &   9066 & $- 9.9$ \\
 7  & 386692 & 411413 & 441849 & $ 78.4$ \\
 7  & 773384 & 822826 & 883698 & $  9.8$ \\
% 7  & \, 773384 \, & \, 822826 \, & \, 883698 \, & $  9.8$ \\
\multicolumn{5}{c}{}\\
\end{tabular}
\qquad
\begin{tabular}{r|c|c|c|r}
$n$ &    $x$ &   $y$ &   $z$ & $r$\ \
\\ \hline
 8  & 209959 & 629874 & 629886 & $-11.6$ \\
 8  & 209945 & 629826 & 629838 & $ 11.6$ \\
 9  &   6817 &  10727 &  10747 & $  5.3$ \\
 9  &  21860 &  25208 &  25903 & $ 24.7$ \\
10  &    280 &    305 &    316 & $137.1$ \\
10  &    560 &    610 &    632 & $ 17.1$ \\
10  &    840 &    915 &    948 & $  5.1$ \\
10  &   7533 &   8834 &   8999 & $  4.4$ \\
12  &   1782 &   1841 &   1922 & $  6.1$ \\
12  &   3987 &   4365 &   4472 & $- 7.1$ \\
12  & 781769 & 852723 & 874456 & $ 10.3$ \\
13  &    666 &    806 &    811 & $  8.3$ \\
13  &   5579 &   8235 &   8239 & $  4.1$ \\
15  & 434437 & 588129 & 588544 & $ 42.9$ \\
%15 & \, 434437 \, & \, 588129 \, & \, 588544 \, & $ 42.9$ \\
16  & 492151 & 741267 & 741333 & $  4.6$ \\
19  &     79 &     85 &     86 & $- 4.7$ \\
19  &    491 &    565 &    567 & $  4.9$ \\
19  &  43329 &  51144 &  51257 & $  5.8$ \\
20  &   4110 &   4693 &   4709 & $  4.3$ \\
\end{tabular}
}

\vspace*{1ex}

All decimal values of~$r$ are rounded to the nearest tenth.
If for some integer $\lambda>1$ we have $r>4\lambda^3$
then $(\lambda x,\lambda y,\lambda z)$ will also appear
in the table provided $\lambda z \leq 10^6$; this happens for
$\lambda=2$ at $n=4,5,7,10$, and for $\lambda=3$ at $n=5,10$.
The first examples for $n=10$ and particularly $n=5$ (where
$13^5+16^5=17^5+12$) are small and striking enough that one
feels they must have been observed already, but I do not
know a reference.  On the other hand, the first two examples
for $n=12$ have been published, and in a most unlikely place:
each appeared in a different episode of the popular animated
cartoon {\em The Simpsons}.  Perhaps the third example for $n=12$,
or an example with $n=7$ or $n=15$, could be used if the cartoon
repeats this theme once more; the relative error $|z^n-y^n-x^n|/z^n$
in each case is between $1$ and $2$ parts in $10^{18}$, as compared
to $3\cdot 10^{-10}$ and $2\cdot 10^{-11}$ for the two four-digit
examples\ldots

Frivolity aside, one is struck by the pair of examples for $n=8$.
The values of $r$ are far from the largest in the table,
but they are almost equal and opposite, and involve nearly equal
triples $(x,y,z)$ for which $z-y$ has the same small value of~$12$.
This suggests that we are dealing with a polynomial family
$(x(t),y(t),z(t))$ specialized at $t=\pm t_0$.  Indeed we quickly find
that these are the cases $t=\pm 3$ of
\be
(32 t^9 + 6 t)^8 + (32 t^8 + 7)^8 = (32 t^9 + 10 t)^8 +
21\cdot 2^{28} t^{40} + O(t^{32}),
\label{id8}
\ee
with $r=t^5/21+O(t^{-3})$.  Thus arbitrarily large values of~$r$
occur, and indeed $z^8-y^8-x^8$ can be as small as $O(z^{40/9})$
rather than the expected $O(z^5)$.  Trying to generalize the
identity (\ref{id8}) further, we soon find that there are
similar families for any exponent $n$ such that $3n(n-2)$
is a square:

\begin{theorem}
Let $n>1$ be a positive integer.  Then there exist polynomials
$x(t),y(t),z(t)\in\Z[t]$ of the form
\be
x(t) = C t^n + D,\quad
y(t) = A t^{n+1} + B t,\quad
z(t) = A t^{n+1} + B' t
\label{xyz.thm}
\ee
with $A\neq0$, $B'\neq B$\/ such that $z^n-y^n-x^n$ is a polynomial
of degree at most $n(n-3)$, if and only if $3n(n-2)$ is a square.
In that case, there exist infinitely many integer triples $(x,y,z)$
with $0<x<y<z$ such that $z^n-y^n-x^n \ll z^{(n^2-3n)/(n+1)}$.
\label{fermat}
\end{theorem}

\begin{proof}
Let $b,b'$ be the distinct rational numbers $B/A,B'/A$.
Expand $z^n-y^n$ at infinity:
\bea
z^n - y^n &=&
n A^n (b'-b) \left(
 t^{n^2} + \frac{n-1}{2} (b'+b) t^{n^2-n}\right.
\label{xn.taylor} \\
 && \ \left. + \frac{(n-1)(n-2)}{6}({b'}^2+b'b+b^2) t^{n^2-2n}
 + O(t^{n^2-3n})
\right).
\nonumber
\eea
For this to be of the form $(Ct^n+D)^n+O(t^{n^2-3n})$ we must have
\be
(n-1) \left(\frac{n-1}{2}(b'+b)\right)^{\!2} =
2n \frac{(n-1)(n-2)}{6}({b'}^2+b'b+b^2).
\label{xn.quadratic}
\ee
The discriminant of this quadratic equation in $b'/b$
is $3n(n-2)$ times a square; thus (\ref{xn.quadratic})
has nonzero rational solutions if and only if $3n(n-2)$
is a square.  Explicitly we find that $b,b'$ are proportional
to $\sqrt{(n^2-2n)/3} \pm 1$.

Conversely, suppose $n^2-2n=3m^2$ for some integer $m$.  Let
\be
z = A(t^{n+1} + c(m+1)t),\qquad
y = A(t^{n+1} + c(m-1)t).
\label{yz.pf}
\ee
Then
\be
z^n-y^n =
2cn A^n \Bigl( t^n + \frac{n-1}{n} c m \Bigr)^n + O(t^{n^2-3n}).
\label{xn.pf}
\ee
To make this $(Ct^n+D)^n+O(t^{n^2-3n})$ with $C,D\in\Z$ we now
need only choose nonzero $c\in\Z$ so that $2cn$ is an $n$-th
power (e.g.\ take $c=(2n)^{n-1}$), and then choose $A$\/ so that
$n|Acm$.  Specializing $t$ to sufficiently large integers in
the resulting $(x(t),y(t),z(t))$ yields infinitely many integer triples
$(x,y,z)$ with $0<x<y<z$ such that $z^n-y^n-x^n \ll z^{(n^2-3n)/(n+1)}$,
as claimed.
\qed
\end{proof}

The smallest $n>3$ such that $n^2-2n=3m^2$ is $n=8$.  There are
infinitely many further examples, starting with $27$, $98$, $363$,
\ldots, and parametrized by a Fermat-Pell equation.  Dropping the
constraint $n>3$ yields the further cases $n=2$ and $n=3$.  For
$n=2$ we again obtain a Pythagorean parametrization, this time
with $x,y,z$ multiplied by~$t$; for $n=3$ we find
\be
(9t^3+1)^3 + (9t^4)^3 - (9t^4+3t)^3 = 1,
\label{xyz3.1}
\ee
one of infinitely many polynomial solutions of $x^3+y^3-z^3=1$.

\subsection{The Fermat cubic}

Our algorithm applies to the Fermat cubic as it does to the Fermat
curves of higher degree, but we treat it separately both because
the heuristic analysis is subtler and because the problem of
finding small values of $|z^3-y^3-x^3|$ has already attracted some
attention.  We noted that in general we expect the smallest values
of $|z^n-y^n-x^n|$ to be comparable with $z^{n-3}$.  For $n=3$,
we have $z^{n-3}=1$, and of course (given this case of Fermat's
Last Theorem) $|z^3-y^3-x^3|$ can be no smaller than~$1$ for
nonzero integers $x,y,z$.  Moreover, $z^3-y^3-x^3$ cannot be
an arbitrary rational multiple of $z^{n-3}$: only the discrete
values $\pm1,\pm2,\ldots$ may arise.  Thus, instead of a
Diophantine inequality $z^n-y^n-x^n\ll z^{n-3}$, we have
a family of Diophantine equations $z^3-y^3-x^3=d$\/ ($d\in\Z$),
and new tools can bear on solving them or, failing that, describing
their distribution of solutions.  These equations have been
investigated by various means since the beginning of the computer
age; see \cite{Guy} for references to work up to about 1980
(some of which dates back to the 1950's), and
\cite{Bremner,CV,HBLR,KTS,PV} for more recent results.
As we shall see, the problem has been approached in several ways,
some of which already improve on direct exhaustion over some $N^2$
values of $(x,y)$.  Still, our new linear approximation method
is better yet, both in heuristic theory --- even though
by factors smaller than our accustomed $N/\LogN$ --- and in practice,
as evidenced by the computation of many new solutions.  Our discussion
here applies with almost no change to other ``diagonal'' cubics, such as
$x^3+y^3+2z^3$ which was also singled out in~\cite[Prob.~D5]{Guy};
but we have not yet implemented a search for small values of
$|x^3+y^3+2z^3|$ beyond what has already been reported in the literature.

For each nonzero $d$, the expected distribution of solutions of
\be
z^3-y^3-x^3=d
\label{eq.d}
\ee
involves not only considerations of size --- i.e.\
of local behavior at the archimedean place of~$\Q$ --- but also on the
behavior of $z^3-y^3-x^3$ at finite primes~$p$: each $p$ contributes
a local factor $f_p(d)$ that is the ratio of the $p$-adic measure of the
$\Z_p$-points of (\ref{eq.d}) to the average of that measure
as $d$\/ ranges over~$\Z_p$.  For instance, if any of those factors
$f_p(d)$ vanishes, there can be no solutions at all.  It is not hard
to see that the only such local constraint is $d\not\equiv\pm4\bmod9$.
For such $d$, the resulting product over~$p$ was investigated by
Heath-Brown~\cite{HB}.  He showed that the product does not converge
absolutely, but can nevertheless be analyzed and approximated
numerically by comparing $f_p(d)$ with the factor at~$p$ of the Euler
product for $\bigl(\zeta_{\Q(\!\root3\of d)}(s)/\zeta(s)\bigr)^3$
at $s=1$, which differs from $f_p(d)$ by a factor of at most
$1+O(p^{-3/2})$.  The product $\prod_p f_p(d)$ is then seen to diverge
to $+\infty$ if $d$\/ is a cube and to converge to a positive limit
when $d$\/ is neither a cube nor congruent to $\pm4\bmod9$.
Heath-Brown thus conjectured in~\cite{HB} that all nonzero integers
$d\not\equiv\pm4\bmod9$ occur as $z^3-y^3-x^3$ infinitely often.
So far this is only known when $d$\/ is either a cube or twice
a cube, thanks to polynomial parametrizations, which the above
heuristics do not try to account for.  We have already exhibited
polynomial solutions for $d=1,2$.  For many $d\not\equiv4\bmod9$\/
which are neither cubes nor twice cubes, not a single solution
is known for $z^3-y^3-x^3=d$.  Heath-Brown observes~\cite{HB}
that this is not surprising, because for many of these~$d$\/
the expected number of solutions with $z\in[N,10^6 N]$ is positive
but smaller than~$1$.  Guy~\cite{Guy} lists the cases with $d<10^3$
which were open as of 1980, and while the list is now shorter the
question of which integers are the sums of three cubes is not yet
settled even in that range.  For instance, the case $d=30$ was open
until~1999, and had been the smallest open case for several decades.

We have noted already that a direct search finds all small
$|z^3-x^3-y^3|$ with $z<N$\/ in time $N^2\LogN$.  There have
been several improvements on this, all obtained by rewriting
the equation (\ref{eq.d}) as
\be
x^3 + d = z^3 - y^3 = (z-y) (z^2 + yz + y^2).
\label{ferm3trick:x}
\ee
Once $x^3+d$\/ is factored, which takes heuristic time
$N^{o(1)}$, all solutions of (\ref{ferm3trick:x})
can be found by trying each factor of $x^3+d$ for $z-y$.
Given the value of~$d$, this takes time only $N^{1+o(1)}$.
In addition to dealing with only one $d$\/ at a time,
this method has the disadvantage that the $N^{o(1)}$
time required to factor $x^3+d$, though subexponential,
is still considerable.  The advantage of this method
is that it finds all solutions with $x\leq N$, while $y,z$
may be considerably larger, of order up to $N^{3/2}$.
Many of the new solutions found in~\cite{KTS} are of this
type, with $y,z$ large but $z-y$ very small. Heath-Brown
observed that, again given~$d$, the factorization of $x^3+d$\/
can be simplified by a precomputation in $\Z[\root3\of d\,]$,
though the complexity of the precomputation depends unpredictably
on~$d$\/ via the arithmetic of the number field $\Q(\root3\of d\,)$;
this approach was implemented in~\cite{HBLR}.
Note that in effect these methods find rational points near the Fermat
cubic that are close to the tangents to the curve at its inflection
points --- the same tangents that demand special care in our algorithm.
A further variation which we suggested in 1996 is to use the
factorization
\be
z^3 - d = x^3 + y^3 = (x+y) (x^2 + xy + y^2)
\label{ferm3trick:z}
\ee
as follows: fix $x+y$, solve for $z\bmod x+y$, and try each
of the resulting values of~$z$.  Here we only find solutions
with $z$, not $x$, bounded by~$N$, but the advantage is that
factoring costs are greatly diminished.  To find all
cube roots of~$d\bmod x+y$ requires factoring $x+y$, a number
of size~$N$\/ rather than $N^3$; and with enough space to
set up a sieve the factorization can be avoided entirely.
In 1999, Eric Pine, Kim Yarbrough, Wayne Tarrant and Michael Beck,
all graduate students at the University of Georgia, took up this
suggestion, choosing $d=30$, and found the first solution:
\be
30 =  2220422932^3 - 283059965^3 - 2218888517^3
\label{sol30}
\ee

We announced our new algorithm in the same 1996 posting to the
{\tt NMBRTHRY} mailing list, together with results of a search for
solutions with $z<10^7$ and $|d|<10^3$.  We did our search in~{\tt gp},
making our computation easy to program (since {\tt gp} already provides
multiprecision arithmetic and lattice reduction) but far from optimally
efficient.  In 1999, unaware of the work of the Georgia group, we asked
Dan J.~Bernstein for an efficient implementation.  He soon wrote a
C program that found all solutions with $z<3\cdot 10^9$ and $|d|<10^4$,
including (\ref{sol30}) and many others.  Several values of~$d$\/
had not been previously represented as the sum of three cubes.
Detailed results and analysis will appear elsewhere.  As usual,
since we are interested in small~$d$, not all $d\ll N$, the
improvement by a factor $N^{1/13}$ should apply here as well
to find all cases of $|z^3-y^3-x^3| \ll z^{3/13}$ with $z<N$, but
we have not attempted to implement such a computation.

\subsection{Miscellaneous examples}

{\bf Trinomial units.}  One sometimes sees in Olympiad-style mathematics
contests the question ``Is $z^{1/3}$ greater or smaller than
$x^{1/3} + y^{1/3}$?'' for some specific positive integers
$x,y,z$.  Of course this is a challenge only when the sign of the
difference $u_3 := z^{1/3} - (x^{1/3} + y^{1/3})$
cannot be determined by inspection.  In some cases the question
 be settled by applying classical inequalities; for instance
if $a>b>0$ then $(a+b)^{1/3} + (a-b)^{1/3} < 2 a^{1/3}$
by convexity of the cube root.  The general solution is to compute
the norm of~$u_3/z^{1/3}$, an algebraic number of degree~$9$ none of
whose other conjugates is real unless $x=y$.  We find that $u_3$
has the same sign as
\be
\N(x,y,z) := (z-y-x)^3 - 27xyz.
\label{crux}
\ee
Moreover, given the size of $x,y,z$, the smaller $\N(x,y,z)$ is,
the nearer $u_3$ will be to~$0$.  In particular, we would like to have
$\N(x,y,z)=\pm1$, which would make the algebraic integer~$u_3$ a unit.
Thus again we seek rational points close to a plane cubic curve,
here $\N(x,y,z)=0$.  This time the curve is rational: by construction,
it is parametrized by $(x:y:z)=(t^3:(1-t)^3:1)$.  It is thus not smooth,
but its only singularity is the isolated point $(x:y:z)=(1:-1:-1)$
(geometrically a node with complex conjugate tangents), which
does not affect our algorithm.  The three rational points of inflection
at $xyz=z-y-x=0$ do affect our algorithm, but fortunately we are not
interested in the points on their tangent lines, since those are the
points with $xyz=0$.  We thus restrict our attention to the portion
of the curve with $x/z,y/z>1/N$, i.e.\ with $t\gg N^{-1/3}$ and
$1-t\gg N^{-1/3}$ in the rational parametrization.  This takes us
far enough from the inflection points that they cause us no difficulty.

The situation is now much the same as for $z^3-y^3-x^3=d$.  We expect
the number of solutions of $\N(x,y,z)=d$\/ of height up to~$N$\/ to be
proportional to $\log N$\/ times a product of local factors $g_p(d)$.
The only local factor that can vanish is $g_3(d)$, which is nonzero
if and only if $9|d$\/ or $d\equiv\pm1\bmod9$.  We henceforth assume
that $d$\/ is in one of these congruence classes.  We can then
check whether $\prod_p g_p$ converges by comparing it with the
$L$-series of the projective cubic surface $\N(x,y,z)=d\,t^3$.
This in turn depends on the Galois structure of the N\'eron-Severi
group of the surface, which can be determined from the action of Galois
on the lines on that cubic surface, as explained in~\cite{Weil}.
We must be careful here because, unlike $x^3+y^3+z^3=d\,t^3$,
the surfaces $\N(x,y,z)=d\,t^3$ are not smooth: each has an $A_2$
singularity at $(x:y:z:t)=(1:-1:-1:0)$.  Thus each has, not $27$
lines as usual, but $15$, of which $6$ go through the singularity;
see~\cite{BW}.  Explicitly, these are the preimages under the
projection to $(x:y:z)$ of the three coordinate axes and the
two tangents to the curve at $(1:-1:-1)$.  We conclude that,
as with (\ref{eq.d}), $\prod_p g_p(d)$ converges unless $d$\/ is
a cube.  So we expect the number of unparametrized solutions
of height $\leq N$\/ to grow as $\log N$, except when $d$\/
is a cube, when it should grow faster, albeit still as a power
of $\log N$\/ --- perhaps $\log^3 N$, by analogy with Manin's
conjecture for cubic surfaces.

Unlike the case of (\ref{eq.d}), we know of no solutions of
$\N(x,y,z)=\pm 1$ in nonconstant polynomials $x,y,z\in\Z[t]$,
other than the trivial ones with $xyz=0$.  Nevertheless we can
find infinitely many nontrivial integer solutions
parametrized by Fermat-Pell equations, and thus show
that the number of solutions of height $\leq N$\/ is $\gg\log N$.
There are several ways to do this.  In 1982 we found a somewhat
complicated route to such a parametrization, obtaining a family
of solutions starting with $\N(16948,31226,186919)=-1$.  The
details may be found in the pages of~\cite{Crux}.  Many years later,
we observed that a simpler approach is to factor $\N(x,y,z)=\pm 1$ as
\be
27 xyz = (z-y-x)^3 \mp 1
= (z-y-x\mp1) \bigl[(z-y-x)^2 \pm (z-y-x) + 1 \bigr].
\label{cruxtrick}
\ee
For each $r\in\Q^*$, we obtain a conic curve $C_r$ by setting
$(z-y-x\mp1)=rx$ in (\ref{cruxtrick}).  This can be
viewed geometrically as follows: the affine surface $\N(x,y,z)=\pm 1$
contains the line $x=(z-y-x\mp1)=0$; thus the intersection of the
surface with any plane $(z-y-x\mp1)=rx$ containing that line
is the union of the line and some residual conic, which is our $C_r$.
Likewise we could start from the line $z=(z-y-x\mp1)=0$
and intersect it with a variable plane $(z-y-x\mp1)=rz$.
For many choices of~$r$, one of these conics is a hyperbola
with infinitely many integral points
parametrized by a Fermat-Pell equation.

In retrospect this approach to $\N(x,y,z)=\pm 1$, in which we fiber
an affine surface by conics that may be regarded as principal
homogeneous spaces for Fermat-Pell equations, seems a remarkable
premonition of our later analysis~\cite{NDE:Euler} of the projective
quartic surface $A^4+B^4+C^4=D^4$ via a fibration by \hbox{genus-$1$}
curves (principal homogeneous spaces for elliptic curves).  In both
cases the approach finds infinitely many solutions but does not readily
lend itself to efficiently finding all solutions of height $\leq N$.
Again a later computation found that the solution that was discovered
first, because it lies on the first fiber that could contain a
solution, is not the one of smallest height.  We used our algorithm
to find all small values of $\N(x,y,z)$ with $0<x,y,z\leq 10^6$.
We found that the smallest solution of $\N(x,y,z)=\pm1$ is
$(14,84,313)$ of norm $+1$, followed by $(6818, 4996, 46879)$,
$(20388, 4881, 86830)$, and $(2742, 32540, 96843)$ each of norm~$-1$,
the known $(16948,31226,186919)$, and $(3408, 182899, 370338)$ of
norm~$+1$, with no further solutions up to $10^6$.  We also found
several primitive solutions of $\N(x,y,z)=\pm8$ and a few sporadic
examples with $d$\/ small but not a cube, which could not have been
obtained at all using the factorization trick; the smallest of these
are
\be
\N(204, 115327, 162434) = 17,
\qquad
\N(650,1425,7899) = 26.
\label{crux:sporadic}
\ee
The $\N(x,y,z)=17$ solution yields a disappointingly large value
of $u_3$ because the conjugates
$z^{1/3} - y^{1/3} - e^{\pm 2\pi i/3} x^{1/3}$
are smaller than usual.  An unexpected result --- since the
identity (\ref{fermdiff}) cannot be used with exponents $<1$
--- was a polynomial solution of $\N(x,y,z)=108$, namely
$(4,y(t),-y(-1-t))$ where $y(t)=4t^3-6t+3$.  We can write
this symmetrically as $(8,g(t),-g(-t))$ where $g(t) = 8t^3-12t-6t+11$,
a cubic polynomial determined up to scaling by the condition that
the Laurent expansion at infinity of $(g(t))^{1/3}$ have
vanishing $t^{-2}$ and $t^{-4}$ terms.  In this form, $\N(x,y,z)$
is the larger constant $864=2^3 108$, but with the bonus that $x$
is a cube so $u_3$ involves one fewer surd; for instance, taking
$t=7$ we find that $\root 3\of {3279}$ is smaller than
$2 + 5\root3\of{17}$ by less than $3.75\cdot 10^{-7}$.  In this
family, as with the first example in~(\ref{crux:sporadic}),
$u_3$ is of order $z^{-2}$, not $z^{-8/3}$, because two of
the conjugates of $u_3$ are $O(1)$.  

A similar investigation of $u_4 := z^{1/4} - (x^{1/4} + y^{1/4})$
was not as productive, perhaps not surprisingly since there
are no arithmetic reasons to expect many nonzero small examples.
For the record, the smallest $z^{11/4} |u_4|$ value found for
$z<10^6$ was $0.365+$ for $(x,y,z)=(241,691,6759)$,
while the smallest $|u_4|$ in that range was $(3.23-) \cdot 10^{-16}$
for $(37792, 36109, 591093)$.

{\bf The $\pi$-th Fermat curve.}  To illustrate our algorithm
also for non-algebraic curves, we chose to apply it to the
Fermat curve of exponent~$\pi$.  Since $\pi$ exceeds $3$,
but only slightly, we expected that $|z^\pi-y^\pi-x^\pi|$ achieves
a global minimum over all $x,y,z$ with $0<x\leq y<z$ but that
the minimum might involve numbers of several digits.  We were
rewarded with the example
\be
2063^\pi + 8093^\pi - 8128^\pi = 0.019369- = 8128^{\pi-3}/(184.75+),
\label{pi8128}
\ee
which seems likely to be the minimum of $|z^\pi-y^\pi-x^\pi|$
over all positive integers $x,y,z$.  At any rate, according
to our computations it is the smallest with $z\leq 10^6$.
The ratio $184.75+$ is also the largest in that range, though
there is also
\be
1198^\pi + 4628^\pi - 4649^\pi = -(0.04949+) = -4649^{\pi-3} / (66.794+).
\label{pi4649}
\ee
It will probably be a long time before the question of the minimality
of (\ref{pi8128}) is settled; a weaker but still intractable conjecture
is that there are only finitely many integer solutions of
$|z^\pi-y^\pi-x^\pi|<1$.

{\bf The Klein quartic.}  All our examples so far were Fermat curves,
even though some had unusual exponents $1/3$, $1/4$, $\pi$.  Probably
the best-known projective plane curve that is not a Fermat curve is
the Klein quartic $K(X,Y,Z)=0$, where
\be
K(X,Y,Z) := X^3 Y + Y^3 Z + Z^3 X.
\label{klein}
\ee
We used our algorithm to search for small values of $K(x,y,z)$.
By symmetry we may assume $\max(x,y,z)=z$.  We are then seeking
rational points near a segment of a plane curve with a single
inflection point, at $x=y=0$.  The tangent $x=0$ at this point
accounts for the obvious family $(0,1,z)$ with \hbox{$K(x,y,z)=z$}.
Our computation up to height~$10^6$ quickly revealed a less obvious
family, $K(1,-t^2,t^3)=-t^2$, with $K(x,y,z)$ growing even more
slowly than the height.  As usual we also found sporadic examples,
though here (as with several other cases we have already seen
such as the Fermat quintic) the best ones are small enough that
our algorithm was not needed to locate them:
\bea
K(1421,-1057,1501) &=& -49, \nonumber \\
K(7211,-8381,11010) &=& -121, \label{klein:sporadic} \\
K(-1550,11817,32615) &=& 245, \nonumber
\eea
with $z/|K(x,y,z)| = 30.6$, $91.0$, $133.1$ respectively.
The largest $z/|K(x,y,z)|$ found with $z\in[10^5,10^6]$
off the singular cubic $y^3+x^2 z = 0$ was $6.756+$, from
$K(-7871, 175577, 829244)=122741$.

\section{Hall's conjecture}

\subsection{Review of Hall's conjecture}

By {\em Hall's conjecture} me mean the following assertion:
if $x,y$ are positive integers such that
\be
k:=x^3-y^2
\label{k.def}
\ee
is nonzero (equivalently, such that $(x,y)\neq(t^2,t^3)$), then
\be
|k| \gg_\eps x^{1/2-\eps}.
\label{Hallconj}
\ee
(While this accords with current usage,
it is not exactly what Hall originally wrote:
as F.~Beukers points out, Hall~\cite{Hall} conjectured
$|k| \gg x^{1/2}$, a stronger statement which is probably false
--- the usual heuristic suggests that there are at least
$(\delta+o(1))\log X$\/ cases of $0<|k|<\delta\sqrt x$ with
$x<X$\/ --- but unlikely to be soon disproved.  See also
\cite{BCHS} for the early history of this conjecture.)
Among several equivalent forms of (\ref{Hallconj}) 
we note the conjecture that the discriminant of
an elliptic curve over~$\Q$ in its standard minimal
form has absolute value $\gg_\eps |a_4|^{1/2-\eps}$.
Known lower bounds on~$|k|$ are much weaker than~(\ref{Hallconj}).
By Siegel's theorem on the finiteness of integer points on elliptic
curves, each nonzero $k\in\Z$ occurs finitely many times as $x^3-y^2$,
so $|k|\ra\infty$ as $x\ra\infty$.  Siegel's proof is ineffective
and thus says nothing about how fast $|k|$ must grow with~$x$.
Starting with Baker's method, effective bounds have become available,
but they are still very weak.  For instance, it is not yet possible
to prove for any $\theta>0$ that $|k|\gg x^\theta$.

Hall's conjecture is now recognized as an important special
case of the Masser-Oesterl\'e ABC conjecture \cite{Oesterle}
(see also \cite{Lang}).  Thus its analogue over function fields
is known to be true by Mason's theorem~\cite{Mason}.  In the special
case of Hall's conjecture for polynomials $x(t),y(t)$, the fact that
$x^3-y^2$ is either zero or has degree $>\frac12\deg(x)$ was proved
some twenty years earlier by Davenport~\cite{Davenport} in response
to a question raised in~\cite{BCHS}.
As in~\cite{NDE:ABC} it follows that the conjecture cannot be
{\em dis}proved by a polynomial parametrization, and indeed in any
polynomial family $(x(t),y(t)|t\in\Z)$ we must have $k\gg x^\theta$
with $\theta>1/2$.  One does better with solutions parametrized
by Fermat-Pell equations, i.e.\ $x,y\in\Z[t,\sqrt{at^2+bt+c}]$
for some $a,b,c\in\Z$ such that $u^2=at^2+bt+c$ has infinitely many
solutions.  The function field $\Q(t,\sqrt{at^2+bt+c})$ is then still
rational, so the Davenport-Mason inequality again holds, but since now
there are two places at infinity one can have $x^3-y^2$ of degree
exactly $\frac12\deg(x)$, and thus attain $\theta=1/2$.  The
existence of a single such family (exhibited below) shows that
the exponent in~(\ref{Hallconj}) cannot be raised above $1/2$.
The fact that one cannot reduce $\theta$ below $1/2$ in this way
was again observed in~\cite{NDE:ABC} in the more general context
of the ABC conjecture.  This fact lends some credence to that
conjecture, and thus to its special case (\ref{Hallconj}); this
contrasts with the situation for $|z^n-y^n-x^n|$, where there is no
reason why some polynomial or Pell family might not do better than the
$z^{n-3-\eps}$ expected by probabilistic heuristics, and indeed
we found such families for some choices of~$n$.

We next digress to say some more on polynomial and Fermat-Pell families
that attain the Davenport-Mason bound, both because they are of
independent interest and because families of both kinds appear
in our numerical results.  In either case $x^3-y^2=k$ is an identity
in a genus-zero function field, namely $\Q(t)$ in the polynomial case
and $\Q(t,\sqrt{at^2+bt+c})$ in the Fermat-Pell case.  Let
$x,y$ have degrees $2m,3m$ respectively, and suppose $k$\/
has the smallest degree possible, i.e.\ $m+1$ in the polynomial
case and $m$ for Fermat-Pell.  Then $f := x^3/y^2$ is a rational
function of degree~$6m$ or~$12m$ on~$\Pr^1$ ramified only above
$0,1,\infty$.  The Riemann existence theorem provides infinitely many
such functions $f=x^3/y^2$ in $\C(t)$; this answers the first part of
the question raised in \cite[p.68]{BCHS}.  The second part concerns
solutions over~$\R$, and can probably be settled by adding data on
complex conjugation to the branched covering.  But we are most interested
in the third part of the question, in which $f$\/ must have rational
coefficients.  Given any one $(x(t),y(t))$, we may trivially obtain
others of the form $(x',y')=(\lambda^2 x(t'),\lambda^3 y(t'))$ where
$t'=at+b$ in the polynomial case, and $t'\in\Q[t]$ with
$\sqrt{a{t'}^2+bt'+c} \big/ \sqrt{at^2+bt+c} \in\Q[t]$
in the Fermat-Pell case.  If we regard such $(x',y')$ and $(x,y)$
as equivalent, only a handful of examples over~$\Q$ are known,
and there may well be no others.  We next list representatives of
the known examples.

In the polynomial case, all known examples have $m\leq 5$.  For $m=1$,
translation and scaling brings any quadratic $x(t)$ to the form
$t^2+2a$, and then $y=t^3+3at$ and $k=3a^2t^2+8a^3$. Necessarily
$a\neq0$, and all such examples are ``twists'' of each other,
becoming isomorphic over $\bar\Q$ if not over~$\Q$.  Note that
$x^3/y^2$ is a degree-$3$ function of $t^2$ with a triple zero.
This function occurs for instance as the cover of the modular curve
$\X(1)$ by $\X_0(2)$.
For $m=2$ we again find that the solution is unique up to twist:
$x=t^4+4at$, $y=t^6+6at^3+6a^2$, and $k=-8a^3t^3-36a^4$.  This time
$x^3/y^2$ is a degree-$4$ function of $t^3$, whose ramification
identifies it with the modular cover $\X_0(3)\ra\X(1)$.
Birch found examples of $(x,y,k)$ with $m=3,5$ and included them in a
29.ix.1961 letter to Chowla; they are reported in \cite{BCHS}:
{
\small
\be
%(x,y,k) =
\left(
36t^6 + 24t^4 + 10t^2 + 1,\
216 t^9 + 216 t^7 + 126 t^5 + 35 t^3 + \frac{21}{4} t,\
\frac92 t^4 + \frac{39}{16} t^2 + 1
\right),
\label{Birch3}
\ee
}%
and
{
\small
\be
%(x,y,k) =
\left(
\frac{t}{9}(t^9 + 6t^6 + 15t^3 + 12),\
\frac{t^{15}}{27} + \frac{t^{12} + 4t^9 + 8t^6}{3} + \frac{5t^3+1}{2},\
-\frac{3t^6 + 14t^3 + 27}{108}
\right).
\label{Birch5}
\ee
}%
These yield integer solutions if $t$ is a multiple of~$4$
in~(\ref{Birch3}) or congruent to~$3$ mod~$6$ in~(\ref{Birch5}).
As noted in~\cite{BCHS}, the second example provides infinitely
many integer solutions of $|x^3-y^2|\ll x^{3/5}$; moreover, for
this choice of twist, the leading coefficient of $k(t)$ is small enough
that $|x^3-y^2|$ is even a respectably small multiple of~$x^{1/2}$
for the first few specializations of~$t$.  The maps $f=x^3/y^2$
associated with Birch's polynomials both have interesting Galois
groups.  For~(\ref{Birch3}), $f$ is a degree-$9$ function of~$t^2$
whose Galois group is $\PSL_2(\F_8)$ over~$\C(t^2)$ and
$\Aut(\PSL_2(\F_8))$ over~$\Q(t^2)$; the Galois closure is the
Fricke-Macbeath curve \cite{Fricke,Macbeath}.  For~(\ref{Birch5}),
$f$ is a degree-$10$ function of~$t^3$ whose Galois group is
$\PSL_2(\F_9)$.  These groups and curves do not arise in connection
with classical modular curves, but they can be identified with certain
Shimura modular curves, most naturally those associated with with
the $(2,3,7)$ and $(2,3,8)$ arithmetic triangle groups (see for
instance~\cite{Takeuchi,NDE:Shimura}).  Hall~\cite[p.185]{Hall}
gives an example with $m=4$:
\be
x = 4(t^8 + 6t^7 + 21t^6 + 50t^5 + 86t^4 + 114t^3 + 109t^2 + 74t + 28);
\label{Hall4}
\ee
In August 1998 I announced a new example with $m=5$ (its computation
will be explained elsewhere):
{
\small
\be
x =
t^{10}+2t^9+33t^8+12t^7+378t^6-336t^5+2862t^4-2652t^3+14397t^2-9922t+18553.
\label{NDE5}
\ee
}%
In both cases (as with all the other $(x,y,k)$ examples), $y$ is
obtained by truncating the Laurent expansion at infinity of $x^{3/2}$
after the constant term.  Neither (\ref{Hall4}) nor (\ref{NDE5}) yields
an interesting Galois group: the Galois groups of $x^3/y^2$ are
${\rm Alt}_{24}$ and $\Sym_{30}$ respectively.  While (\ref{NDE5}),
like (\ref{Birch5}), must yield infinitely many integer solutions of
$|x^3-y^2|\ll x^{3/5}$, the leading coefficient of~$k$ in (\ref{NDE5})
makes the implied constant much larger, and none of these solutions
will appear in our list of small values of $|x^3-y^2|$.  The question,
raised in~\cite{BCHS}, whether there are any $x,y,k\in\Q[t]$ of degrees
$2m,3m,m+1$ with $m>5$, remains unsolved.

For Fermat-Pell families, the list is even shorter: all known
examples are equivalent, and come from the identity
\be
(t^2+10t+5)^3 - (t^2+22t+125) (t^2+4t-1)^2 = 1728t.
\label{X5}
\ee
Here $y$ is a multiple of $\sqrt{at^2+bt+c}$, so $f$\/ factors
as a map of degree~$6$ composed with the double cover of~$\Q(t)$
by $\Q(t,\sqrt{at^2+bt+c})$.  We noted in \cite[p.49]{NDE:Modular}
that the resulting degree-$6$ map $f=x^3/y^2: \Pr^1\ra\Pr^1$
is the cover $\X_0(5)\ra\X(1)$ of classical modular curves.
Thus the elliptic curves of low discriminant coming from
the identity~(\ref{X5}) all admit a rational $5$-isogeny.
Each Fermat-Pell family obtained from~(\ref{X5}) by specifying
the class of $t^2+22t+125$ mod ${\Q^*}^2$ yields
$k\sim Cx^{1/2}$ for some nonzero~$C$.  The smallest such $C$\/
is $5^{-5/2} 54 = .96598\ldots$, obtained by Danilov~\cite{Danilov}
by substituting $125(2t-1)$ for~$t$ in~(\ref{X5}) and dividing by $20^3$:
\be
(5^5 t^2 - 3000 t + 719)^3
- (5^3 t^2 - 114 t + 26) (5^6 t^2 - 5^3 123 t + 3781)^2 = 27(2t-1).
\label{Danilov}
\ee
The factor $5^3 t^2 - 114 t + 26$ is a square for $t=-5$, and thus
for infinitely many~$t$.  The first case $t=-5$ of this yields
the elliptic curve of discriminant~$-11$ labeled 11-A2(C) in
Cremona's table~\cite{Cremona}; it is known that the isogeny
class of this curve provides the examples with minimal conductor
of a rational $5$-isogeny, and indeed of an elliptic curve over~$\Q$.

\subsection{The new algorithm}

To obtain numerical data with which to compare Hall's conjecture,
we want to find all small nonzero values of $|x^3-y^2|$ with
$x,y\in\Z$ and $x\leq X$.  So that we can compare our algorithm
with other approaches we briefly review previous work on this problem.

The most direct approach is to simply compute for each $x\leq X$\/
the integer~$y$ closest to $x^{3/2}$.  Since $x^{3/2}$ varies
smoothly with~$x$, this can be done quite efficiently, but clearly
must take at least time proportional to~$X$.  This is essentially
what Hall did in~\cite{Hall}, with $X=7\cdot 10^8$; some three
decades later, faster computers make larger~$X$\/ feasible, and
indeed Frits Beukers reports in an Aug.~1998 e-mail that he performed
such a computation for $X=10^{12}$.  But this is probably close to
the practical limit with today's technology, and at any rate this
direct approach is superseded by the $X^{1/2}\LogX$\/ algorithm
described below.

A fundamentally different approach is taken in~\cite{GPZ}:
for each nonzero $k\in[-K,K]$, investigate the arithmetic
of the elliptic curve $E_k:y^2=x^3-k$, and use effective bounds
on integral points to find all integer solutions of $x^3-y^2=k$.
In~\cite{GPZ}, Gebel, Peth\"o, and Zimmer did most of this work
for~$K=10^5$, except for a few values of~$k$, for which they could not
find a generator for $E_k(\Q)$; Wildanger later showed in his doctoral
thesis~\cite{Wildanger} that none of these $E_k$ has an integral point,
thus completing the computation of integer solutions of $0<|x^3-y^2|<10^5$.
It is not clear even heuristically how this method compares with
other approaches.  It is the only approach used thus far that will
provably find all solutions with $|k|\leq K$.  (The recent proof
of the modularity conjecture means that Cremona's
algorithms~\cite{Cremona} yield another such approach,
but to my knowledge it has not been used to solve $x^3-y^2=k$.)
Assuming Hall's conjecture, $|k|\leq K$\/ is equivalent to
$x \ll_\eps K^{2+\eps}$, but this begs the question
of the constant implied in~``$\ll$''$\!$.  Neither do we know how to
estimate the average work required to find all integer points on
a curve $E_k$.  It may be reasonable to guess that this average work
is proportional to $K^{c+o(1)}$ for some $c>0$.  (This estimate
certainly holds for Cremona's algorithms.)  The total work would then be
$K^{1+c+o(1)}$.   Under Hall's conjecture, this is equivalent to
$X^{(1+c)/2+o(1)}$, so strictly worse (modulo an unknown implied
constant) than our $X^{1/2}\LogX$\/ algorithm, though perhaps better
than a direct search, depending on whether $c<1$.

We noted already that the direct search can exploit the smoothness
of the function $x^{3/2}$.  We can try to take further advantage
of this by mimicking our approach to rational approximation of curves:
surround the segment $x<X$\/ of the semicubical parabola $y=x^{3/2}$
by a union of parallelograms each of area~$O(1)$, and use lattice
reduction to quickly find all integer points in each parallelogram.
This does give an asymptotic improvement, though a small one: the
parallelogram containing a point $(x,x^{3/2})$ has length $\gg x^{1/6}$,
so the computational cost is reduced by at most $X^{1/6}$, to
$X^{5/6}\LogX$.

We reduce the exponent of~$X$\/ from~$1$ or~$5/6$ to $1/2$ by a
more radical reorganization of the computation that lets us apply
lattice reduction more efficiently.  More generally, for each
positive $c\in\Q$\/ we can find all cases of $0<|cx^3-y^2|\ll x$\/
in time $O_c(X^{1/2}\LogX)$.  All choices of~$c$ are essentially
equivalent: we get from one to the other by scaling $x,y$ and imposing
congruence conditions on them.  The most convenient choice of~$c$ turns
out to be $4/3$.  We thus show how to solve $0<|4x^3-3y^2|\ll x$;
the cases relevant to Hall's conjecture are those with $3|x$ and $6|y$,
when $(4x^3-3y^2)/108=(x/3)^3-(y/6)^2$.

We begin as in~\cite{Hall} by approximating~$x,y$ by (multiples of)
a square and a cube.  Any positive integer $x$ may be written
uniquely as
\be
x = 3 \zeta^2 + \eta
\quad{\rm with}\quad
\eta,\zeta\in\Z,
\zeta>0,\,\eta\in(-3\zeta,3\zeta].
\label{z.eta}
\ee
Then
\be
(4x^3/3)^{1/2} = 6 \zeta^3 + 3 \eta \zeta + \frac14 \frac{\eta^2}{\zeta}
- \frac1{72} \left(\frac\eta\zeta\right)^{\!3} + O(1/\zeta).
\label{y:taylor}
\ee
We thus write
\be
y = 6 \zeta^3 + 3\eta \zeta + \xi
\label{xi}
\ee
with $\xi\ll\zeta$.  More precisely, if
\be
\eta = \beta \zeta
\label{beta}
\ee
Then $\beta\in(-3,3]$, and $|4x^3-3y^2|\ll x$ if and only if
\be
\xi = \frac{\eta^2}{4\zeta} - \frac1{72} \beta^3 + O(1/\zeta).
\label{xi:taylor}
\ee
At this point, Hall~\cite{Hall} imposes the assumption
$\beta\ll\xi^{-1/5}$.  We allow an arbitrary $\beta\in(-3,3]$
and approximate it within $O(X^{-1/2})$ by one of $O(X^{1/2})$
evenly spaced points in that interval.  Suppose, then, that 
$b$ is one of those points.  We approximate (\ref{xi:taylor})
by a linear combination of $\zeta$, $\eta-b\zeta$, and~$1$:
\be
\xi = \frac{b^2}{4} \zeta + \frac{b}{2} (\eta-b\zeta) - \frac{b^3}{72}
+ O(1/\zeta)
= - \frac{b^2}{4} \zeta + \frac{b}{2} \eta - \frac{b^3}{72}
+ O(1/\zeta).
\label{xi:linear}
\ee
We now assume that $\zeta\gg X^{1/2}$, for instance by requiring that
$x>X/4$; repeating the computation with $X$\/ replaced by $X/4$, $X/16$,
$X/64$,\ldots\ will then cover the entire range $x\leq X$, and if we
can cover $(X/4,X]$ in time $O(x^{1/2}\LogX)$ then the same is true of
$[1,X]$.  Under the assumption $x\in(X/4,X]$,
we have the following constraints on $\xi,\eta,\zeta$:
\be
\zeta \ll X^{1/2}, \quad
\eta-b\zeta \ll 1,
\label{z.eta.cond}
\ee
and
\be
\xi + \frac{b^2}{4} \zeta - \frac{b}{2} \eta + \frac{b^3}{72}
\ll X^{-1/2}.
\label{xi.cond}
\ee
We are thus in a familiar situation: we seek all the integral points
in $O(X^{1/2})$ parallelepipeds, each of volume $O(1)$.  The term
$b^3/72$ in~(\ref{xi.cond}) means that the parallelepipeds are
no longer centered at the origin, but this causes no difficulty ---
indeed we already dealt with off-center parallelepipeds in the
practical implementation of our algorithm for finding rational points
near curves.  So again we linearly transform each parallelepiped to a
cube and obtain a lattice reduction problem; if these lattices were
randomly distributed among three-dimensional lattices, we would almost
certainly have only $O(X^{1/2})$ points to try, and would thus find
all solutions of $0<|4x^3-3y^2|\ll x$ with $x\leq X$\/ in time
$O(X^{1/2}\LogX)$.

In fact it turns out that in this case our lattices are {\em not}\/
equidistributed: they all lie in a $2$-dimensional subspace of the
$5$-dimensional moduli space of lattices in~$\R^3$.  This gives rise
to both a minor annoyance and a major advantage.  The bad news is
that we cannot expect our lattices to have on average $O(1)$ vectors
of norm $\ll 1$; but this annoyance is minor because the actual
average is proportional to $\log X$\/ and thus can be absorbed
into the $\LogX$\/ factor.  The good news is that we understand
our special lattices well enough to actually prove results
that are only heuristic for rational points near curves.

The key is that in each case our lattice is a {\em symmetric square}\/
of a lattice in~$\R^2$.  By this we mean the following.  Recall that
the symmetric square of a \hbox{$2$-dimensional} vector space~$V$\/
is the $3$-dimensional vector space $\Sym^2 V$\/ consisting of
symmetric tensors in $V\otimes V$.  Since $\Sym^2 V$\/ is defined
naturally in terms of~$V$, any linear transformation of~$V$\/
yields a linear transformation of~$\Sym^2 V$.  We thus have a
homomorphism $\Sym^2: \GL_2\ra\GL_3$.  To give this map explicitly
we choose a basis $(e_1,e_2)$ for~$V$, and use the basis
$(e_1\otimes e_1, (e_1\otimes e_2 + e_2\otimes e_1)/2, e_2\otimes e_2)$
for~$\Sym^2 V$.  We then calculate that
\be
\Sym^2 \left(\begin{array}{cc}p&q\\r&s\end{array}\right) =
\left(
\begin{array}{ccc}
p^2 &   p q   & q^2 \\
2pr &\;ps+qr\;& 2qs \\
r^2 &   r s   & s^2
\end{array}
\right) .
\label{sym2.def}
\ee
Over any field, $\Sym^2(\SL_2)$ is contained in the subgroup of $\SL_3$
preserving the discriminant form $4a_1 a_3 - a_2^2$ on $\Sym^2(V)$; if
we worked over an algebraically closed field, that subgroup would
coincide with $\Sym^2(\SL_2)$.  Now (\ref{z.eta.cond},\ref{xi.cond})
mean that the column vector $v=(\xi,\eta,\zeta)\in\Z^3$ satisfies
$\|M_b v - u_b\| \ll 1$ where $u_b := (0,0,-b^3/72)$ and
\be
M_b := \left(
\begin{array}{ccc}
   0    &          0          & X^{-1/2} \\
   0    &          1          &   -b     \\
X^{1/2} &\;\, -X^{1/2}b/2 \,\;& X^{1/2}b^2/4
\end{array}
\right)
= \Sym^2 \left(\begin{array}{cc}
0 & X^{-1/4} \\ X^{1/4} & X^{1/4} b/2
\end{array}
\right).
\label{mb.def}
\ee
This is why we went after $4x^3-3y^2$ rather than pursuing $x^3-y^2$
directly: an analogous approach to $x^3-y^2$ would yield a matrix
that is still a symmetric square but with respect to a different
basis, requiring a definition of $\Sym^2$ with fractional coefficients
and complicating the lattice reduction.  Note that the quadratic form
$4\xi\zeta-\eta^2$ preserved by $M_b\in\Sym^2(\SL_2)$ is already
visible in~(\ref{xi:taylor}).

Our algorithm, then, is as follows.  For each of our $O(X^{1/2})$
choices of~$b$, calculate the matrix
\be
N_b := \left(\begin{array}{cc}
0 & X^{-1/4} \\ X^{1/4} & X^{1/4} b/2
\end{array}
\right)
\label{nb.def}
\ee
with $M_b = \Sym^2 = N_b$.  Use lattice reduction to find a matrix
$K_b\in\GL_2(\Z)$ such that $N_b K_b$ is as small as possible.
Then
\be
M'_b := \Sym^2(N_b K_b) = \Sym^2 N_b \; \Sym^2 K_b = M_b \; \Sym^2 K_b.
\label{m'b.def}
\ee
is small too.  Let $L_b = \Sym^2 K_b \in \GL_3(\Z)$.
Then $M_b v = M'_b L_b^{-1} v$.  Find a box containing all
$w\in\Z^3$ such that $\|M'_b w - u_b\| \ll 1$.  For each $w$ in
the box, compute $v=L_b w$ and check whether the resulting
$x,y$ satisfy $x\in(X/4,X]$ and $0<|4x^3-3y^2|\ll x$; if they do,
output $x$ (and check whether $3|x$ and $6|y$ to determine whether
this solution also yields a small value of $x^3-y^2$).  This is
easier than our usual algorithm because we are reducing a lattice
in~$\R^2$ rather than $\R^3$, which in our case amounts to calculating
the continued fraction of $b/X^{1/2}$.  Moreover, the computational
cost of the algorithm can be bounded rigorously: $M'_b$ will only
be large if $b/X^{1/2}$ is close to a rational number with numerator
and denominator $\ll X^{1/4}$, and the effect of such a close
rational approximation is easy to determine.  Summing over all
rationals of height $\ll X^{1/4}$ we find that the total number
of candidate vectors~$v$ is $\ll X^{1/2}\log X$, and thus that
the computation takes time $O(X^{1/2}\LogX)$ as claimed.

Note that the $X^{1/2}\log X$\/ bound also has the following
consequence: there are $\ll X^{1/2} \log X$\/ solutions of
$|x^3-y^2| \ll x$\/ with $x\leq X$.  Moreover, if $C$\/
is large enough, we can deduce from this analysis that there
are $\gg X^{1/2}$ solutions of $0<|x^3-y^2|<Cx$ with $x\in[X/2,X]$.
More generally, we show that for each positive $c\in\Q$ there exists
$C$\/ such that for each $r\in\R/\Z$ and $d>1$ there are at most
$C d X^{1/2} \log X$ solutions of $|(cx^3)^{1/2}-(y+r)| < d X^{-1/2}$
with $x,y\in\Z$ and $x<X$\/; and, given $c$ as above and any
$\theta\in[0,1)$, there exists $C_0$ such that for any $r\in\R/\Z$
there are $\gg X^{1/2}$ solutions of $|(cx^3)^{1/2}-(y+r)|<C_0 X^{-1/2}$
with $x,y\in\Z$ and $x\in[\theta X,X]$.  The constants $C,C_0$
depend effectively on $c,\theta$.  These results improve considerably
on results in this direction available from general exponential-sum
techniques for proving uniform distribution mod~$1$.  The detailed
proofs of our claims in this paragraph will appear elsewhere.

\subsection{Numerical results}

We have implemented our algorithm in a C program using $64$-bit
integer arithmetic, again replacing each $O(\cdots)$ and $\ll$ by
explicit bounds, and searched for all solutions of
$0<|4x^3-3y^2|<200x^{1/2}$ with $4\cdot 10^6 < x < 3\cdot 10^{18}$.
The range $x<10^{10}$ was covered by a direct search, the overlap
$[4\cdot 10^6,10^{10}]$ being used as a check on the computation.
The code was processed with an optimizing compiler and ran for
three weeks during the summer of~1998 on a Sun Sparcstation Ultra~$1$.
As a corollary we obtained all cases of $0<|x^3-y^2|<\frac12\sqrt x$
with $x<10^{18}$.  (With currently available hardware the same
computation could easily finish in a few days; with parallelization
it should be feasible to reach $10^{23}$ at least.)
The next table lists, for each of the $25$ solutions of
$0<|x^3-y^2|<\frac12\sqrt x$, the values of $k=x^3-y^2$, $x$, 
and $r=x^{1/2}/|k|$.  We need not list $y$, which is always
the integer nearest to $x^{3/2}$.  The explanation of the last
two columns follows the table.

\vspace*{1ex}

\centerline{
\begin{tabular}{c|r|r|r|c|c}
\# & $k$\quad{} & $x$\qquad{} & $r$\quad{} & GPZ? & Comments
\\ \hline
 1 & 1641843 & 5853886516781223 & 46.60 &   & !! \\
 2 & 30032270 & 38115991067861271 & 6.50 &   & ! \\
 3 & $-$1090 & 28187351 & 4.87 & + & \\
 4 & $-$193234265 & 810574762403977064 & 4.66 &   & \\
 5 & $-$17 & 5234 & 4.26 & + & $P(-3)$ \\
 6 & $-$225 & 20114 & 3.77 & + & \\
 7 & $-$24 & 8158 & 3.76 & + & $P(3)$ \\
 8 & 307 & 939787 & 3.16 & + & \\
 9 & 207 & 367806 & 2.93 & + & \\
10 & $-$28024 & 3790689201 & 2.20 & + & \\
11 & $-$117073 & 65589428378 & 2.19 &   & \\
12 & $-$4401169 & 53197086958290 & 1.66 &   & \\
13 & 105077952 & 23415546067124892 & 1.46 &   & * \\
14 & $-1$ & 2 & 1.41 &   & \\
15 & $-$497218657 & 471477085999389882 & 1.38 &   & \\
16 & $-$14668 & 384242766 & 1.34 & + & $P(-9)$ \\
17 & $-$14857 & 390620082 & 1.33 & + & $P(9)$ \\
18 & $-$87002345 & 12813608766102806 & 1.30 &   & \\
19 & 2767769 & 12438517260105 & 1.27 &   & \\
20 & $-$8569 & 110781386 & 1.23 & + & \\
21 & 5190544 & 35495694227489 & 1.15 &   & \\
22 & $-$11492 & 154319269 & 1.08 & + & \\
23 & $-$618 & 421351 & 1.05 & + & \\
24 & 548147655 & 322001299796379844 & 1.04 &   & D \\
25 & $-$297 & 93844 & 1.03 & + & D
\end{tabular}
}

\vspace*{1ex}

The ``GPZ'' column indicates whether the solution was among
the $13$ listed in~\cite{GPZ}.  These are the solutions
with $1<|k|<10^5$.  Presumably the solution $2^3-3^2=-1$
is not on that list because the elliptic curve $y^2=x^3+1$
was already known to have rank~$0$ so Gebel, Peth\"o and Zimmer
were not interested in it.

The \#1 row is a new record, improving the previous record $r$
by a factor of almost~$10$, whence the notation ``!!''.  Even
row~\#2, marked ``!'', has $r$ larger than the old record
which is row~\#3.  Either of this suffices to refute
Hall's comment \cite[p.175]{Hall}, repeated in~\cite{GPZ},
that $r<5$ seems to hold in all cases.

*: Obtained from row \#1 by scaling $(x,y,k)$ to $(2^2 x,2^3 y,2^6 k)$.
This reduces $r$ by a factor of~$32$, but $r=46+$ in row~\#1
is large enough that even $r/32$ still exceeds the threshold
of our table.

$P(t)$: Birch's polynomial family (\ref{Birch5}).  This has
$r=12/t+O(t^{-4})$, so the only values of $t\equiv 3\bmod 6$
that appear on the $r>1$ list are $t=\pm3$ and $\pm9$.
Already in~\cite[p.69]{BCHS} the specializations $t=\pm3$
are noted as ``striking special cases'' of~(\ref{Birch5}).

D: The first two cases of Danilov's family~(\ref{Danilov}).  The
appearance of the larger of these was a welcome check on our
computation.

\vspace*{1ex}

Any threshold on~$r$ is of necessity arbitrary; the next solution
has $r$ just below our cutoff of~$1$:
$(x,k,r)=(16544006443618,4090263,0.9944\ldots)$.

\end{document}